# Transform the non-linear programming problem to the initial-value problem to solve

Sheng ZHANG, Fei LIAO, Yi-Nan KONG, and Kai-Feng HE

**Abstract** A dynamic method to solve the Non-linear Programming (NLP) problem with Equality Constraints (ECs) and Inequality Constraints (IECs) is proposed. Inspired by the Lyapunov continuous-time dynamics stability theory in the control field, the optimal solution is analogized to the stable equilibrium point of a finite-dimensional dynamic system and it is solved in an asymptotic manner. Under the premise that the Karush-Kuhn-Tucker (KKT) optimality condition exists, the Dynamic Optimization Equation (DOE), which has the same dimension to that of the optimization parameter vector, is established and its solution will converge to the optimal solution of the NLP globally with a theoretical guarantee. Using the matrix pseudo-inverse, the DOE is valid even without the linearly independent regularity requirement on the nonlinear constraints. In addition, the analytic expressions of the Lagrange multipliers and KKT multipliers, which adjoin the ECs and the IECs respectively during the entire optimization process, are also derived. Via the proposed method, the NLP may be transformed to the Initial-value Problem (IVP) to be solved, with mature Ordinary Differential Equation (ODE) integration methods. Illustrative examples are solved and it is shown that the dynamic method developed may produce the right numerical solutions with high efficiency.

**Keywords** Non-linear programming; Karush-Kuhn-Tucker optimality condition; Lyapunov dynamics stability; Dynamic optimization equation; Initial-value problem

**Subject classification codes** 90C30

## 1. Introduction

Non-linear Programing (NLP) aims to determine the parameters that optimize a specified objective function while satisfying various Equality Constraints (ECs) and Inequality Constraints (IECs). In the physical world, many scientific and engineering problems may be abstracted as NLP problems, and the general formulation is defined as

***Problem 1*:** Consider the objective function

$$J = f(\boldsymbol{\theta}), \tag{1}$$

subject to

$$\boldsymbol{g}(\boldsymbol{\theta}) \leq \boldsymbol{0}, \tag{2}$$

The authors are with the Computational Aerodynamics Institution, China Aerodynamics Research and Development Center, Mianyang, 621000, China. (e-mail: zszhangshengzs@hotmail.com).

$$\boldsymbol{h}(\boldsymbol{\theta}) = \boldsymbol{0}, \tag{3}$$

where $\boldsymbol{\theta} \in \mathbb{R}^n$ is the optimization parameter vector. $f : \mathbb{R}^n \to \mathbb{R}$ is a scalar function with continuous first-order partial derivatives with respect to $\boldsymbol{\theta}$. $\boldsymbol{g} : \mathbb{R}^n \to \mathbb{R}^r$ and $\boldsymbol{h} : \mathbb{R}^n \to \mathbb{R}^s$ are $r$-dimensional vector function and $s$-dimensional vector function with continuous first-order partial derivatives, respectively. Find the optimal solution $\hat{\boldsymbol{\theta}}$ that minimizes $J$.

Theories and methods on the computation of the NLP have been widely studied [1, 2]. Because of their complexity, generally NLPs are solved with numerical methods. Traditional methods usually use the iterative mechanism to seek the solution. Prevailing methods include the Sequential Quadratic Programming (SQP) method and the Interior-point (IP) method. The SQP method is considered to be one of the most efficient methods for the constrained optimization. At every iteration, an approximate Quadratic Programming (QP) sub-problem is solved, and the solution is gradually achieved with a sequence of QP sub-problems [3]. Since its proposal in Wilson [4], this method has been systematically developed to achieve the desired convergence for optimization with general nonlinear constraints [5-8]. The IP method is also a popular method and it employs the barrier parameters to treat the IECs. The resulting sub-problems, corresponding to the decreasing barrier parameters, are solved and their solutions converge to the solution of the original problem [9]. There has been a better understanding of the IP methods [10, 11] and efficient algorithms have been developed with desirable performance [12-14].

Besides the numerical iterative methods, there is another way to solve the NLP, which is based on the continuous-time dynamics. With such method, the Dynamic Optimization Equation (DOE), which is a set of differential equation, is developed and the optimization problem is transformed to the Initial-value Problem (IVP) to be solved. Studies on this type of methods may date back to the 1940s for the unconstrained problems [15]. Other relevant work regarding the unconstrained minimization includes the gradient dynamic equation (e.g. Refs. [16-19]), the second-order dynamic equation arising from the physical energy view [20], and the continuous Newton method [21], etc. To address the constrained NLP problems, Brown and Bartholomew-Biggs [22] utilize the penalty function to facilitate the application of the dynamic method for unconstrained problems. Tanabe [23] establishes the DOE for the NLP with ECs in the feasible solution region, and it may also address the IECs by transforming them to the ECs with the quadratic slack parameters. Yamashita [24] further generalizes the DOE for the NLP with ECs in the infeasible solution region. Evtushenko and Zhadan [25] also present the similar DOE; moreover, with the slack parameter and coordinate transformation, NLP with both ECs and IECs is reformulated as NLP with ECs only to be solved. By introducing the quadratic slack parameters to equalize the IECs, Schroop [26] develops the corresponding DOE and an equivalent index-2 Differential Algebraic Equation (DAE) form is also presented.

Despite its convenience, employing quadratic slack parameters to circumvent the tough IECs not only increases the parameters to be determined, but may also greatly increase the number of minimums compared with the original problem [27]. To address such issue, Jongen and Stein [28] first introduce the gradient vector field with respect to the original optimization parameter for the NLP with IECs only, in the light of Riemannian metric upon the work of Rapcsák [29]. Shikhman and Stein [30] further derive the DOE defined in the original optimization parameters for the NLP with both ECs and IECs in the feasible solution region, through computing the projected gradient in the original, lower-dimensional space from the higher-dimensional space with slack parameters. DOE proposed under the framework of the IP method may also avoid the employment of slack parameters. For example, based on the transformed merit function, Moguerza and Prieto [31] combine three search directions to get the approximate solution of the DOE system upon the gradient flow; Ali and Oliphant [32] establish the DOE, which is motivated by the work of Snyman [20], to solve the joint primal-dual problem. However, such methods introduce extra dynamics of the multiplier parameters as the cost.

Recently, inspired by the Lyapunov continuous-time dynamics stability theory in the control field, a dynamic method, the Variation Evolving Method (VEM), is proposed to solve the Optimal Control Problems (OCPs) [33-35]. Since the NLP problem may be considered as the static case of the OCP, here the dynamic method for the NLP defined in Problem 1 is developed under the similar framework. Compared with the former work, our main contributions are: i) The DOEs for the NLPs are established from the unified view of the Lyapunov principle. As an alternative realization to Shikhman and Stein [30], the IECs are treated directly upon their dynamic attribute in this study. No slack parameters are introduced and the resulting DOE has a dimensionality same to that of the optimization parameter vector. ii) The equation is further generalized to be effective in the infeasible solution region, and more importantly, its solution will converge to the optimal solution of the NLP globally with a theoretical guarantee. iii) Expressions for the Lagrange multipliers and the Karush-Kuhn-Tucker (KKT) multipliers at non-stationary points or even infeasible points are derived, and they will converge to the right solutions as the optimization parameters tend to the optimal value. iv) In addition, most of the aforementioned articles have a linearly independent regularity assumption for the nonlinear constraints, while here the matrix pseudo-inverse is used to address the singularity arising from constraints dependence and the right solution may still be sought without the Linearly Independent Constraint Qualification (LICQ).

Throughout the paper, our work is built upon the assumption that the solution for the NLP problem exists. Study regarding the existence of solutions is beyond the scope of this paper. In the following, first preliminaries including the optimality condition, the matrix pseudo-inverse, the Lyapunov stability theory, and the motivation for the proposed method are presented in Section 2. In Section 3, the DOE that seeks the optimal solution of the NLP within the feasible solution region is derived. In Section 4, the equation is modified to be effective in the infeasible solution region through the first-order stable dynamics principle. A Lyapunov function is constructed to ensure its validity theoretically and strategies handling unsolvable situation are presented. Later in Section 5, illustrative examples are solved to verify the effectiveness of the method. Section 6 concludes the paper with some final remarks.

## 2. Preliminaries

### *2.1. Optimality condition*

We first give necessary definitions for the further study.

**Definition 2.1:** The feasible solution region $\mathbb{D}_f$ is the collection of parameters that satisfy Eqs. (2) and (3), i.e., $\mathbb{D}_f = \{\boldsymbol{\theta} \mid \boldsymbol{g}(\boldsymbol{\theta}) \le \boldsymbol{0}, \boldsymbol{h}(\boldsymbol{\theta}) = \boldsymbol{0}\}$.

With Definition 1, Problem 1 may be expressed concisely as

$$\hat{\boldsymbol{\theta}} = \underset{\boldsymbol{\theta} \in \mathbb{D}_f}{\arg\min}(J) . \quad (4)$$

**Definition 2.2:** The infeasible solution region $\mathbb{D}_{if}$ is the collection of parameters that violate Eqs. (2) or (3).

**Definition 2.3:** For the *i*-th IEC in Eq. (2), it is said to be an optimal active IEC for Problem 1 if $g_i(\hat{\boldsymbol{\theta}}) = 0$.

Note that the optimal active IEC is defined on the optimal solution, while other IECs may still be activated or even violated during the optimization process. The KKT optimality condition characterizes the optimal solution and provides the basis for algorithms. In practice, its existence is ensured by the Constraint Qualification (CQ) in that the linear approximation of the constraints at the optimal solution adequately captures the geometry of the feasible set [1]. The following assumption is made first in our study.

**Assumption 2.1:** The CQ holds in Problem 1 so that the KKT optimality condition exists.

Upon Assumption 2.1, the Lagrange function for Problem 1 may be presented with the multiplier technique, i.e.,

$$\bar{J} = f(\boldsymbol{\theta}) + {\boldsymbol{\pi}_E}^{\mathrm{T}} \boldsymbol{h} + {\boldsymbol{\pi}_I}^{\mathrm{T}} \boldsymbol{g} , \quad (5)$$

where $\boldsymbol{\pi}_E$ is the Lagrange multiplier parameter vector and $\boldsymbol{\pi}_I$ is the KKT multiplier parameter vector. The superscript "$\mathrm{T}$" denotes the transpose operator. Then the first-order differential relation may be derived as

$$\mathrm{d}\bar{J} = (f_{\boldsymbol{\theta}} + {\boldsymbol{g}_{\boldsymbol{\theta}}}^{\mathrm{T}} \boldsymbol{\pi}_I + {\boldsymbol{h}_{\boldsymbol{\theta}}}^{\mathrm{T}} \boldsymbol{\pi}_E)^{\mathrm{T}} \mathrm{d}\boldsymbol{\theta} + \boldsymbol{h}^{\mathrm{T}} \mathrm{d}\boldsymbol{\pi}_E + \boldsymbol{g}^{\mathrm{T}} \mathrm{d}\boldsymbol{\pi}_I , \quad (6)$$

where $f_{\boldsymbol{\theta}} = \dfrac{\partial f}{\partial \boldsymbol{\theta}}$ is the shorthand notation in the form of column vector. $\boldsymbol{g}_{\boldsymbol{\theta}}$ and $\boldsymbol{h}_{\boldsymbol{\theta}}$ are the Jacobi matrixes. Through $\mathrm{d}\bar{J} = 0$ and the analysis on the property of IECs, we may get the conditions that determine the stationary points, including the feasibility conditions (2), (3) and the KKT optimality condition, namely

$$\begin{array}{ll} f_{\theta} + \boldsymbol{h}_{\theta}{}^{\mathrm{T}}\boldsymbol{\pi}_E + \boldsymbol{g}_{\theta}{}^{\mathrm{T}}\boldsymbol{\pi}_I = \boldsymbol{0} & \\ (\pi_I)_i \geq 0 \quad i \in \hat{\mathbb{I}} & , \\ (\pi_I)_i = 0 \quad i \notin \hat{\mathbb{I}} & \end{array} \tag{7}$$

where $\hat{\mathbb{I}}$ is the index set of the optimal active IECs for Problem 1, which is defined as

$$\hat{\mathbb{I}} = \{ i \mid g_i(\hat{\boldsymbol{\theta}}) = 0,\ i = 1,2,...,r \} . \tag{8}$$

In particular, regarding the optimality condition (7), $\boldsymbol{\pi}_E$ and $\boldsymbol{\pi}_I$ will be unique if the LICQ holds, i.e., the gradients of the EC and the optimal active IEC functions are linearly independent. Otherwise, there may be multiple solutions for the multipliers. Note that in the paper we do not assume the LICQ, which is a more conservative requirement [36], and allow the multiple solutions case for the multipliers.

### *2.2. Pseudo-inverse of matrix*

The matrix pseudo-inverse [37] is a generalization of the matrix inverse, and it may handle non-invertible cases.

**Definition 2.4:** The Moore-Penrose pseudo-inverse of a matrix $\boldsymbol{M}$ may be defined as

$$\boldsymbol{M}^{+} = \boldsymbol{V}\boldsymbol{\Sigma}^{+}\boldsymbol{U}^{\mathrm{T}} , \tag{9}$$

where $\boldsymbol{M} = \boldsymbol{U}\boldsymbol{\Sigma}\boldsymbol{V}^{\mathrm{T}}$ is the Singular Value Decomposition (SVD) of $\boldsymbol{M}$. $\boldsymbol{U}$ and $\boldsymbol{V}$ are orthogonal matrixes. $\boldsymbol{\Sigma} = \begin{bmatrix} \boldsymbol{\Sigma}_r & \boldsymbol{0}_{\Sigma 12} \\ \boldsymbol{0}_{\Sigma 21} & \boldsymbol{0}_{\Sigma 22} \end{bmatrix}$ and $\boldsymbol{\Sigma}^{+} = \begin{bmatrix} \boldsymbol{\Sigma}_r^{-1} & \boldsymbol{0}_{\Sigma 21}{}^{\mathrm{T}} \\ \boldsymbol{0}_{\Sigma 12}{}^{\mathrm{T}} & \boldsymbol{0}_{\Sigma 22}{}^{\mathrm{T}} \end{bmatrix}$. $\boldsymbol{\Sigma}_r$ is a diagonal positive-definite block matrix. $\boldsymbol{0}_{\Sigma 12}$, $\boldsymbol{0}_{\Sigma 21}$, and $\boldsymbol{0}_{\Sigma 22}$ are the right-dimensional zero block matrixes.

Obviously, when $\boldsymbol{M}$ is invertible, there is $\boldsymbol{M}^{+} = \boldsymbol{M}^{-1}$, and for any matrix $\boldsymbol{M}$

$$\boldsymbol{M}^{+} = (\boldsymbol{M}^{\mathrm{T}}\boldsymbol{M})^{+}\boldsymbol{M}^{\mathrm{T}} = \boldsymbol{M}^{\mathrm{T}}(\boldsymbol{M}\boldsymbol{M}^{\mathrm{T}})^{+} . \tag{10}$$

By investigating the geometric meaning of SVD, the projection matrix $\boldsymbol{P}_{\boldsymbol{M}c}$, corresponding to the column space $\mathbb{S}_{\boldsymbol{M}c}$ spanned by the column vectors of $\boldsymbol{M}$, is $\boldsymbol{P}_{\boldsymbol{M}c} = \boldsymbol{M}\boldsymbol{M}^{+}$; the projection matrix $\boldsymbol{P}_{\boldsymbol{M}r}$, corresponding to the row space $\mathbb{S}_{\boldsymbol{M}r}$ spanned by the row vectors of $\boldsymbol{M}$, is $\boldsymbol{P}_{\boldsymbol{M}r} = \boldsymbol{M}^{+}\boldsymbol{M}$. Moreover, $\boldsymbol{1} - \boldsymbol{M}\boldsymbol{M}^{+}$ and $\boldsymbol{1} - \boldsymbol{M}^{+}\boldsymbol{M}$ are also projection matrixes, corresponding to the spaces orthogonal to $\mathbb{S}_{\boldsymbol{M}c}$ and $\mathbb{S}_{\boldsymbol{M}r}$, respectively. Here $\boldsymbol{1}$ denotes the right-dimensional identity matrix.

**Remark 2.1:** For two matrixes $\boldsymbol{M}_1$ and $\boldsymbol{M}_2$, $\boldsymbol{P}_{\boldsymbol{M}_1 c} = \boldsymbol{P}_{\boldsymbol{M}_2 c}$ when they have same column space, and $\boldsymbol{P}_{\boldsymbol{M}_1 r} = \boldsymbol{P}_{\boldsymbol{M}_2 r}$ when they have same row space.

Upon the properties of the projection matrix, now the pseudo-inverse in solving linear equations is introduced.

**Lemma 2.1** [37]**:** *For the linear equation*

$$\boldsymbol{M}\boldsymbol{x}=\boldsymbol{b}\,, \tag{11}$$

*where the matrix* $\boldsymbol{M}$ *is arbitrary,*

$$\boldsymbol{M}^{+}\boldsymbol{b}=\underset{\boldsymbol{x}\in\arg\min(\|\boldsymbol{M}\boldsymbol{x}-\boldsymbol{b}\|_2)}{\arg\min}(\|\boldsymbol{x}\|_2)\,, \tag{12}$$

*which means that* $\boldsymbol{M}^{+}\boldsymbol{b}$ *is the optimal solution for the objective* $\|\boldsymbol{x}\|_2$ *within the set* $\left\{\boldsymbol{x}\mid\boldsymbol{x}\in\mathbb{R}^m,\boldsymbol{x}=\arg\min(\|\boldsymbol{M}\boldsymbol{x}-\boldsymbol{b}\|_2)\right\}$*. Here* $\|\cdot\|_2$ *denotes the 2-norm of vector.*

By considering the variants of the linear equation (11), i.e.,

$$\boldsymbol{S}\boldsymbol{M}\boldsymbol{x}=\boldsymbol{S}\boldsymbol{b}\,, \tag{13}$$

and

$$\boldsymbol{M}\boldsymbol{S}(\boldsymbol{S}^{-1}\boldsymbol{x})=\boldsymbol{b}\,, \tag{14}$$

where $\boldsymbol{S}$ is the right-dimensional nonsingular square matrix, we may further derive

**Remark 2.2:** Consider the linear equation (11) where the matrix $\boldsymbol{M}$ is arbitrary. There is

$$(\boldsymbol{S}\boldsymbol{M})^{+}\boldsymbol{S}\boldsymbol{b}=(\boldsymbol{M}^{\mathrm{T}}\boldsymbol{Q}\boldsymbol{M})^{+}\boldsymbol{M}^{\mathrm{T}}\boldsymbol{Q}\boldsymbol{b}=\underset{\boldsymbol{x}\in\arg\min(\|\boldsymbol{S}(\boldsymbol{M}\boldsymbol{x}-\boldsymbol{b})\|_2)}{\arg\min}(\|\boldsymbol{x}\|_2)\,, \tag{15}$$

where $\boldsymbol{Q}=\boldsymbol{S}^{\mathrm{T}}\boldsymbol{S}$ is a positive-definite matrix. Particularly when exact solutions for Eq. (11) exist, there is

$$\boldsymbol{M}^{+}\boldsymbol{b}=(\boldsymbol{M}^{\mathrm{T}}\boldsymbol{M})^{+}\boldsymbol{M}^{\mathrm{T}}\boldsymbol{b}=(\boldsymbol{M}^{\mathrm{T}}\boldsymbol{Q}\boldsymbol{M})^{+}\boldsymbol{M}^{\mathrm{T}}\boldsymbol{Q}\boldsymbol{b}=(\boldsymbol{S}\boldsymbol{M})^{+}\boldsymbol{S}\boldsymbol{b}\,. \tag{16}$$

**Remark 2.3:** Consider the linear equation (11) where the matrix $\boldsymbol{M}$ is arbitrary. There is

$$\boldsymbol{S}(\boldsymbol{M}\boldsymbol{S})^{+}\boldsymbol{b}=\boldsymbol{Q}\boldsymbol{M}^{\mathrm{T}}(\boldsymbol{M}\boldsymbol{Q}\boldsymbol{M}^{\mathrm{T}})^{+}\boldsymbol{b}=\underset{\boldsymbol{x}\in\arg\min(\|\boldsymbol{M}\boldsymbol{x}-\boldsymbol{b}\|_2)}{\arg\min}(\|\boldsymbol{S}^{-1}\boldsymbol{x}\|_2)\,, \tag{17}$$

where $\boldsymbol{Q}=\boldsymbol{S}^{\mathrm{T}}\boldsymbol{S}$ is a positive-definite matrix.

### *2.3. Lyapunov stability theory*

The Lyapunov stability theory investigates the dynamic behaviour of states within a dynamic system, from the view of generalized energy [38].

**Definition 2.5:** For a continuous-time autonomous dynamic system like

$$\dot{\boldsymbol{x}}=\boldsymbol{f}(\boldsymbol{x})\,, \tag{18}$$

where $\boldsymbol{x} \in \mathbb{D} \in \mathbb{R}^n$ is the state, $\dot{\boldsymbol{x}} = \frac{\mathrm{d}\boldsymbol{x}}{\mathrm{d}t}$ is its time derivative, and $\boldsymbol{f} : \mathbb{D} \rightarrow \mathbb{R}^n$ is a vector function. $\mathbb{D}$ is a certain set. If $\hat{\boldsymbol{x}} \in \mathbb{D}$ satisfies $\boldsymbol{f}(\hat{\boldsymbol{x}}) = \boldsymbol{0}$, then $\hat{\boldsymbol{x}}$ is called an equilibrium point.

**Definition 2.6:** The equilibrium point $\hat{\boldsymbol{x}}$ is an asymptotically stable equilibrium point in $\mathbb{D}$, if for any initial condition $\boldsymbol{x}(t)\big|_{t=0} = \boldsymbol{x}_0 \in \mathbb{D}$, there is $\lim\limits_{t \to +\infty} \left\| \boldsymbol{x}(t) - \hat{\boldsymbol{x}} \right\|_2 = 0$.

**Lemma 2.2:** (see Khalil [38], with small adaptation) *For the continuous-time autonomous dynamic system* (18)*, if there exists a continuously differentiable (If not, only except at* $\hat{\boldsymbol{x}}$ *) function* $V : \mathbb{D} \rightarrow \mathbb{R}$ *such that*

*i)* $V(\hat{\boldsymbol{x}}) = c$ *and* $V(\boldsymbol{x}) > c$ *in* $\mathbb{D} / \{\hat{\boldsymbol{x}}\}$,

*ii)* $\dot{V}(\boldsymbol{x}) \leq 0$ *in* $\mathbb{D}$ *and* $\dot{V}(\boldsymbol{x}) < 0$ *in* $\mathbb{D} / \{\hat{\boldsymbol{x}}\}$.

*where* $c$ *is a constant. Then* $\boldsymbol{x} = \hat{\boldsymbol{x}}$ *is an asymptotically stable equilibrium point in* $\mathbb{D}$.

For example, maybe $\boldsymbol{f}(\boldsymbol{x})$ in the dynamic system (18) satisfies $(\boldsymbol{x} - \hat{\boldsymbol{x}})^{\mathrm{T}} \boldsymbol{f}(\boldsymbol{x}) < 0$ for any $\boldsymbol{x} \neq \hat{\boldsymbol{x}}$, and then a feasible Lyapunov function can be constructed as

$$V = \frac{1}{2} (\boldsymbol{x} - \hat{\boldsymbol{x}})^{\mathrm{T}} (\boldsymbol{x} - \hat{\boldsymbol{x}}). \tag{19}$$

The dynamics given by $\boldsymbol{f}(\boldsymbol{x})$ determines that $\dot{V} \leq 0$ and $\boldsymbol{x}$ will converge to the equilibrium $\hat{\boldsymbol{x}}$. Figure 1 sketches the trajectory of some state in the stable dynamic system and the corresponding Lyapunov function value. No matter what the initial condition $\boldsymbol{x}_0$ is, as long as it falls into the stable region $\mathbb{D}$, the state $\boldsymbol{x}$ will approach the equilibrium $\hat{\boldsymbol{x}}$ gradually. Meanwhile, the "energy" of the dynamic system, measured by the function $V$, will reach its minimum. Note that the constant $c$ in Lemma 2.2 is allowed to take non-zero value and this facilitates the analogy with the objective function, whose value may not vanish at the optimum.

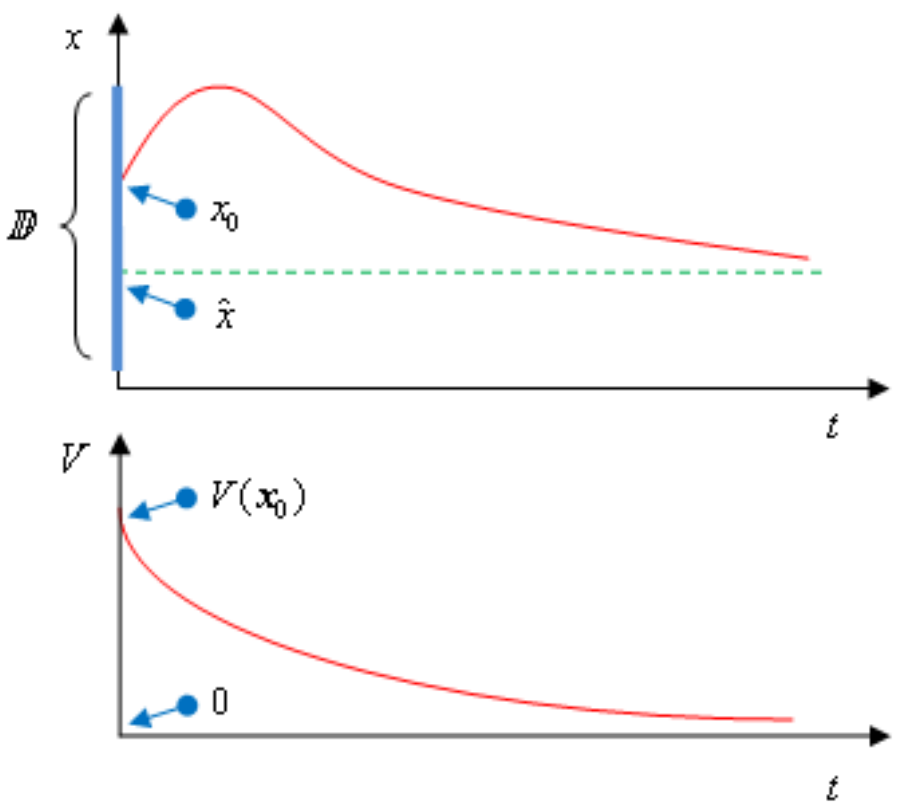

Figure 1. Sketch for the state trajectory and the Lyapunov function value profile.

### *2.4. Motivation*

In the system dynamics theory, from the stable dynamics of state $\boldsymbol{x}$, we may construct a monotonically decreasing "energy" function $V(\boldsymbol{x})$, which achieves its minimum when $\boldsymbol{x}$ reaches $\hat{\boldsymbol{x}}$. Inspired by it, now we consider its inverse problem, that is, from an objective function (regarded as the measure of generalized energy) to derive the dynamics that minimize this objective.

Consider the unconstrained version of Problem 1 where $\boldsymbol{g}$ and $\boldsymbol{h}$ vanish. To find the optimal value $\hat{\boldsymbol{\theta}}$ that minimizes $J$, we make an analogy with the Lyapunov function and introduce a virtual time $\tau$, which is used to describe the derived dynamics. Differentiating $J$ with respect to $\tau$ gives

$$\frac{\mathrm{d}J}{\mathrm{d}\tau} = \frac{\mathrm{d}f}{\mathrm{d}\tau} = {f_{\boldsymbol{\theta}}}^{\mathrm{T}} \frac{\mathrm{d}\boldsymbol{\theta}}{\mathrm{d}\tau}. \tag{20}$$

To guarantee that $J$ decreases with respect to $\tau$, i.e., $\frac{\mathrm{d}J}{\mathrm{d}\tau} \le 0$, we may set the following DOE as

$$\frac{\mathrm{d}\boldsymbol{\theta}}{\mathrm{d}\tau} = -\boldsymbol{K}_{\boldsymbol{\theta}} f_{\boldsymbol{\theta}}, \tag{21}$$

where $\boldsymbol{K}_{\boldsymbol{\theta}}$ is an $n \times n$ dimensional positive-definite matrix. According to Lemma 2.2, under this dynamics, $J$ will decrease until it reaches an minimum, and $\boldsymbol{\theta}$ will approach $\hat{\boldsymbol{\theta}}$ that satisfies $f_{\boldsymbol{\theta}}\big|_{\boldsymbol{\theta}=\hat{\boldsymbol{\theta}}} = \boldsymbol{0}$, the first-order optimality condition. To get the numerical solution of $\hat{\boldsymbol{\theta}}$, the mature Ordinary Differential Equation (ODE) numerical integration methods may be employed to solve Eq. (21) with some initial condition of $\boldsymbol{\theta}\big|_{\tau=0}$. Note that although only the first-order optimality conditions are explicitly satisfied, the solution will not halt at a maximum or a saddle point (unless they are the initial guesses of $\boldsymbol{\theta}$), because those points are not stable equilibrium points. Regarding the virtual time $\tau$, even if it does not explicitly exist, it may be better accepted when referring to the numerical computation dimension. Namely, the numerical iterations performed in the traditional iterative methods are the reflection of the discrete dynamics along the virtual time.

Along this idea, when further considering the general NLP defined in Problem 1, the problem of how to address the IECs (2) and the ECs (3) arises, and the solution will be detailed in the following.

## 3. DOE in feasible solution region

In order to simplify the study, we restrict our consideration in the feasible solution region $\mathbb{D}_f$ first. Under this premise, Eq. (1) will act as the Lyapunov function in deriving the DOE.

### *3.1. NLP with ECs only*

To start with, we consider the NLP with ECs only. Again, to seek the dynamics that ensure the achievement of the optimal solution, Eq. (1) is differentiated with respect to the virtual time $\tau$ to give Eq. (20). However, we cannot use Eq. (21) now because in that way the ECs (3) are not guaranteed. When considered in $\mathbb{D}_f$, we need to find the dynamics that not only guarantees $\frac{\mathrm{d}J}{\mathrm{d}\tau} \leq 0$ but also satisfies the differential relation of

$$\frac{\mathrm{d}\boldsymbol{h}}{\mathrm{d}\tau} = \boldsymbol{h}_{\boldsymbol{\theta}} \frac{\mathrm{d}\boldsymbol{\theta}}{\mathrm{d}\tau} = \boldsymbol{0} , \tag{22}$$

which maintains the feasibility. More importantly, we need to confirm that the solution under such dynamics may reach the optimal solution of the NLP. The answer is given by the following theorem, which is first presented from the analysis of manifold in Tanabe [23]. Here it will be proved upon the Lyapunov principle.

**Theorem 3.1:** *For the NLP with objective* (1) *and ECs* (3), *assume that* $\boldsymbol{h}_{\boldsymbol{\theta}}$ *is row full-rank; solve the IVP defined by the following DOE*

$$\frac{\mathrm{d}\boldsymbol{\theta}}{\mathrm{d}\tau} = -\boldsymbol{K}_{\boldsymbol{\theta}}(f_{\boldsymbol{\theta}} + \boldsymbol{h}_{\boldsymbol{\theta}}{}^{\mathrm{T}} \boldsymbol{\pi}_E) , \tag{23}$$

*with arbitrary initial condition* $\boldsymbol{\theta}\big|_{\tau=0} \in \mathbb{D}_f$, *where* $\boldsymbol{\pi}_E \in \mathbb{R}^s$ *is computed by*

$$\boldsymbol{\pi}_E = -(\boldsymbol{h}_{\boldsymbol{\theta}} \boldsymbol{K}_{\boldsymbol{\theta}} \boldsymbol{h}_{\boldsymbol{\theta}}{}^{\mathrm{T}})^{-1} \boldsymbol{h}_{\boldsymbol{\theta}} \boldsymbol{K}_{\boldsymbol{\theta}} f_{\boldsymbol{\theta}} . \tag{24}$$

$\boldsymbol{K}_{\boldsymbol{\theta}}$ *is an* $n \times n$ *dimensional positive-definite matrix. When* $\tau \to +\infty$, $\boldsymbol{\theta}$ *will satisfy the optimality condition*

$$f_{\boldsymbol{\theta}} + \boldsymbol{h}_{\boldsymbol{\theta}}{}^{\mathrm{T}} \boldsymbol{\pi}_E = \boldsymbol{0} , \tag{25}$$

*and then* $\boldsymbol{\pi}_E$ *also takes the value of* $-(\boldsymbol{h}_{\boldsymbol{\theta}}{}^{\mathrm{T}})^{+} f_{\boldsymbol{\theta}}$, *which is independent of* $\boldsymbol{K}_{\boldsymbol{\theta}}$.

***Proof:*** First, we show that Eq. (23) guarantees that $\frac{\mathrm{d}J}{\mathrm{d}\tau} \leq 0$ and the solution stays in $\mathbb{D}_f$. Though the optimization theory, reformulate Eq. (20) as a constrained optimization problem subject to Eq. (22) to be

$$\begin{gathered} \min \;\; J_{t1} = f_{\boldsymbol{\theta}}{}^{\mathrm{T}} \frac{\mathrm{d}\boldsymbol{\theta}}{\mathrm{d}\tau} \\ \text{subject to} \\ \boldsymbol{h}_{\boldsymbol{\theta}} \frac{\mathrm{d}\boldsymbol{\theta}}{\mathrm{d}\tau} = \boldsymbol{0} \end{gathered} \; . \tag{26}$$

Note that now $\frac{\mathrm{d}\boldsymbol{\theta}}{\mathrm{d}\tau}$ are the optimization parameters. Since the minimum of this problem may be negative infinity, to penalize too large parameters, we introduce another objective function $J_{t2}$ to formulate a Multi-objective Optimization Problem (MOP) as

$$\begin{gathered} \min(J_{t1}, J_{t2}) \\ \text{subject to} \\ \boldsymbol{h}_{\boldsymbol{\theta}} \frac{\mathrm{d}\boldsymbol{\theta}}{\mathrm{d}\tau} = \boldsymbol{0} \end{gathered} , \tag{27}$$

where

$$J_{t2} = \frac{1}{2}\left(\frac{\mathrm{d}\boldsymbol{\theta}}{\mathrm{d}\tau}\right)^{\mathrm{T}} \boldsymbol{K}_{\boldsymbol{\theta}}^{-1} \frac{\mathrm{d}\boldsymbol{\theta}}{\mathrm{d}\tau} . \tag{28}$$

Using the weighting method to solve the Pareto optimal solution of this MOP, the resulting objective is

$$J_{t3} = aJ_{t1} + bJ_{t2} , \tag{29}$$

where $a \geq 0$, $b \geq 0$ and $a + b = 1$. When $a = 1,\ b = 0$, we get a solution that minimizes $J_{t1}$. When $a = 0,\ b = 1$, we get a solution that minimizes $J_{t2}$. Otherwise, we get a compromise solution. For this MOP, obviously in the case of $a = 0,\ b = 1$, the Pareto optimal solution is $\frac{\mathrm{d}\boldsymbol{\theta}}{\mathrm{d}\tau} = \boldsymbol{0}$, and now the values of the objective functions are $J_{t1} = 0$ and $J_{t2} = 0$. For any other cases, the compromise solution guarantees that $J_{t1} \leq 0$ (property of the Pareto optimal solution, see Zitzler [39], Burachik et al. [40] and references therein). Set $a = \frac{1}{2}$ and $b = \frac{1}{2}$; we get the Feasibility-preserving Dynamics Optimization Problem (FPDOP), which is a QP problem as

$$\begin{gathered} \min \quad J_{t3} = \frac{1}{2} J_{t1} + \frac{1}{2} J_{t2} \\ \text{subject to} \\ \boldsymbol{h}_{\boldsymbol{\theta}} \frac{\mathrm{d}\boldsymbol{\theta}}{\mathrm{d}\tau} = \boldsymbol{0} \end{gathered} . \tag{30}$$

Introduce the Lagrange multiplier $\bar{\boldsymbol{\pi}}_E \in \mathbb{R}^r$ to adjoin the constraint. Then the Lagrange function is

$$J_{t4} = \frac{1}{2} J_{t1} + \frac{1}{2} J_{t2} + \frac{1}{2} \bar{\boldsymbol{\pi}}_E^{\mathrm{T}} \boldsymbol{h}_{\boldsymbol{\theta}}^{\mathrm{T}} \frac{\mathrm{d}\boldsymbol{\theta}}{\mathrm{d}\tau} . \tag{31}$$

Let $\frac{\partial J_{t4}}{\partial\left(\frac{\mathrm{d}\boldsymbol{\theta}}{\mathrm{d}\tau}\right)} = \boldsymbol{0}$; we may get the minimum solution for the FPDOP as

$$\frac{\mathrm{d}\boldsymbol{\theta}}{\mathrm{d}\tau} = -\boldsymbol{K}_{\boldsymbol{\theta}} (f_{\boldsymbol{\theta}} + \boldsymbol{h}_{\boldsymbol{\theta}}^{\mathrm{T}} \bar{\boldsymbol{\pi}}_E) . \tag{32}$$

Substitute Eq. (32) into Eq. (22), we have

$$-\boldsymbol{h}_{\boldsymbol{\theta}}\boldsymbol{K}_{\boldsymbol{\theta}}(f_{\boldsymbol{\theta}}+\boldsymbol{h}_{\boldsymbol{\theta}}{}^{\mathrm{T}}\bar{\boldsymbol{\pi}}_{E})=\boldsymbol{0}\ . \tag{33}$$

Since the invertibility of the matrix $\boldsymbol{h}_{\boldsymbol{\theta}}\boldsymbol{K}_{\boldsymbol{\theta}}\boldsymbol{h}_{\boldsymbol{\theta}}{}^{\mathrm{T}}$ is guaranteed by the assumption that $\boldsymbol{h}_{\boldsymbol{\theta}}$ has full row rank, then

$$\bar{\boldsymbol{\pi}}_{E}=-(\boldsymbol{h}_{\boldsymbol{\theta}}\boldsymbol{K}_{\boldsymbol{\theta}}\boldsymbol{h}_{\boldsymbol{\theta}}{}^{\mathrm{T}})^{-1}\boldsymbol{h}_{\boldsymbol{\theta}}\boldsymbol{K}_{\boldsymbol{\theta}}\boldsymbol{f}_{\boldsymbol{\theta}}\ . \tag{34}$$

Furthermore, Eq. (20) may be reformulated as

$$\frac{\mathrm{d}J}{\mathrm{d}\tau}=(f_{\boldsymbol{\theta}}+\boldsymbol{h}_{\boldsymbol{\theta}}{}^{\mathrm{T}}\bar{\boldsymbol{\pi}}_{E})^{\mathrm{T}}\frac{\mathrm{d}\boldsymbol{\theta}}{\mathrm{d}\tau}-\bar{\boldsymbol{\pi}}_{E}{}^{\mathrm{T}}\boldsymbol{h}_{\boldsymbol{\theta}}\frac{\mathrm{d}\boldsymbol{\theta}}{\mathrm{d}\tau}\ . \tag{35}$$

With Eq. (32), and because now Eq. (22) holds, then

$$\frac{\mathrm{d}J}{\mathrm{d}\tau}=-(f_{\boldsymbol{\theta}}+\boldsymbol{h}_{\boldsymbol{\theta}}{}^{\mathrm{T}}\bar{\boldsymbol{\pi}}_{E})^{\mathrm{T}}\boldsymbol{K}_{\boldsymbol{\theta}}(f_{\boldsymbol{\theta}}+\boldsymbol{h}_{\boldsymbol{\theta}}{}^{\mathrm{T}}\bar{\boldsymbol{\pi}}_{E})\ . \tag{36}$$

This means $\frac{\mathrm{d}J}{\mathrm{d}\tau}\leq 0$, and $\frac{\mathrm{d}J}{\mathrm{d}\tau}=0$ occurs only when $f_{\boldsymbol{\theta}}+\boldsymbol{h}_{\boldsymbol{\theta}}{}^{\mathrm{T}}\bar{\boldsymbol{\pi}}_{E}=\boldsymbol{0}$. Besides, we will show $\bar{\boldsymbol{\pi}}_{E}$ introduced here and $\boldsymbol{\pi}_{E}$ used to adjoin the ECs (3) in Eq. (5) are identical at the optimal solution. Compare Eq. (35) with Eq. (6), to achieve the optimality condition within $\mathbb{D}_{f}$, there should be

$$f_{\boldsymbol{\theta}}+\boldsymbol{h}_{\boldsymbol{\theta}}{}^{\mathrm{T}}\bar{\boldsymbol{\pi}}_{E}=f_{\boldsymbol{\theta}}+\boldsymbol{h}_{\boldsymbol{\theta}}{}^{\mathrm{T}}\boldsymbol{\pi}_{E}\ . \tag{37}$$

Since Eq. (37) holds for arbitrary $\boldsymbol{h}$ that has full row rank of $\boldsymbol{h}_{\boldsymbol{\theta}}$, we have $\bar{\boldsymbol{\pi}}_{E}=\boldsymbol{\pi}_{E}$ at the optimal solution.

Now it is easy to show that when $\tau\rightarrow+\infty$, $\boldsymbol{\theta}$ will satisfy the optimality condition (25). By Lemma 2.2, for the dynamic equation (23), Eq. (1) is a Lyapunov function within the feasible solution region $\mathbb{D}_{f}$. Thus, from any initial condition of $\boldsymbol{\theta}\big|_{\tau=0}\in\mathbb{D}_{f}$, the minimum solution that satisfies (25) is an asymptotically stable equilibrium point in $\mathbb{D}_{f}$. Further by Remark 2.2, when Eq. (25) is met, it may be found that $\boldsymbol{\pi}_{E}=-(\boldsymbol{h}_{\boldsymbol{\theta}}{}^{\mathrm{T}})^{+}f_{\boldsymbol{\theta}}=-(\boldsymbol{h}_{\boldsymbol{\theta}}\boldsymbol{h}_{\boldsymbol{\theta}}{}^{\mathrm{T}})^{-1}\boldsymbol{h}_{\boldsymbol{\theta}}f_{\boldsymbol{\theta}}$ for the optimal solution, which is independent of $\boldsymbol{K}_{\boldsymbol{\theta}}$. □

The matrix $\boldsymbol{K}_{\boldsymbol{\theta}}$ in Theorem 3.1 (and in the following) is only required to be positive-definite and it may be either constant or dependent on $\boldsymbol{\theta}$. A constant case is simple to use, while a parameter-dependent one may bring better performance. However, the tuning will be an artful task. The Lagrange multiplier defined by Eq. (24) agrees with the first-order multiplier estimate at a non-stationary point in [3], with $\boldsymbol{K}_{\boldsymbol{\theta}}$ being an identity matrix. However, here it is not an estimate but has definite meaning, and the proof of Theorem 3.1 shows its convergence to the right solution.

**Remark 3.1:** Solve the NLP with objective (1) and ECs (3) in the approach suggested in Theorem 3.1; the Lagrange multiplier given by Eq. (24) will converge to the right solution, i.e., $-(\boldsymbol{h_\theta}^{\mathrm{T}})^{+} f_{\boldsymbol{\theta}}\big|_{\boldsymbol{\theta}=\hat{\boldsymbol{\theta}}}$ as $\boldsymbol{\theta}$ converges to $\hat{\boldsymbol{\theta}}$.

The DOE (23) requires a matrix inverse of $(\boldsymbol{h_\theta} \boldsymbol{K_\theta} \boldsymbol{h_\theta}^{\mathrm{T}})^{-1}$. In practice, it may occur that $\boldsymbol{h_\theta}$ violates the row full-rank assumption even if the constraints in $\boldsymbol{h}$ are mutually independent. Then Eq. (24) will not be applicable. To address such problem, we employ the pseudo-inverse of matrix. With it, the row full-rank assumption in Theorem 3.1 may be removed and it could be modified as

**Theorem 3.2:** *For the NLP with objective* (1) *and ECs* (3)*, under the premise of Assumption 2.1, solve the IVP defined by the DOE* (23) *with arbitrary initial condition* $\boldsymbol{\theta}\big|_{\tau=0} \in \mathbb{D}_f$*, where* $\boldsymbol{\pi}_E \in \mathbb{R}^s$ *is computed by*

$$\boldsymbol{\pi}_E = -(\boldsymbol{h_\theta} \boldsymbol{K_\theta} \boldsymbol{h_\theta}^{\mathrm{T}})^{+} \boldsymbol{h_\theta} \boldsymbol{K_\theta} \boldsymbol{f_\theta} . \tag{38}$$

$\boldsymbol{K_\theta}$ *is an* $n \times n$ *dimensional positive-definite matrix. When* $\tau \to +\infty$*,* $\boldsymbol{\theta}$ *will satisfy the optimality condition* (25)*, and then* $\boldsymbol{\pi}_E$ *also takes the value of* $-(\boldsymbol{h_\theta}^{\mathrm{T}})^{+} f_{\boldsymbol{\theta}}$*, which is independent of* $\boldsymbol{K_\theta}$*.*

***Proof:*** We only need to show that with $\boldsymbol{\pi}_E$ computed by Eq. (38), the results of Eqs. (23) and (25) are the same to those calculated upon $\underline{\boldsymbol{h}}_{\boldsymbol{\theta}}$, in which the redundant linearly dependent rows are removed from $\boldsymbol{h_\theta}$. Our statement is true if

$$\boldsymbol{h_\theta}^{\mathrm{T}} (\boldsymbol{h_\theta} \boldsymbol{K_\theta} \boldsymbol{h_\theta}^{\mathrm{T}})^{+} \boldsymbol{h_\theta} = \underline{\boldsymbol{h}}_{\boldsymbol{\theta}}^{\mathrm{T}} (\underline{\boldsymbol{h}}_{\boldsymbol{\theta}} \boldsymbol{K_\theta} \underline{\boldsymbol{h}}_{\boldsymbol{\theta}}^{\mathrm{T}})^{-1} \underline{\boldsymbol{h}}_{\boldsymbol{\theta}} . \tag{39}$$

By Remark 2.1, there is

$$\boldsymbol{h_\theta}^{\mathrm{T}} (\boldsymbol{h_\theta} \boldsymbol{h_\theta}^{\mathrm{T}})^{+} \boldsymbol{h_\theta} = \underline{\boldsymbol{h}}_{\boldsymbol{\theta}}^{\mathrm{T}} (\underline{\boldsymbol{h}}_{\boldsymbol{\theta}} \underline{\boldsymbol{h}}_{\boldsymbol{\theta}}^{\mathrm{T}})^{-1} \underline{\boldsymbol{h}}_{\boldsymbol{\theta}} . \tag{40}$$

Furthermore, by considering $\boldsymbol{h_\theta} (\boldsymbol{K_\theta})^{1/2}$ and $\underline{\boldsymbol{h}}_{\boldsymbol{\theta}} (\boldsymbol{K_\theta})^{1/2}$ together, we may prove the validity of Eq. (39). Similarly by Remark 2.2, when the optimality condition (25) is met, the argument that $\boldsymbol{\pi}_E = -(\boldsymbol{h_\theta}^{\mathrm{T}})^{+} f_{\boldsymbol{\theta}}$ at the optimal solution may be proved. □

From the optimization theory, we know that when $\boldsymbol{h_\theta}$ does not have full row rank, the solutions for the multiplier are multiple. Yet with the value given by Eq. (38), $\boldsymbol{\pi}_E$ is optimal in the sense of minimum 2-norm at $\hat{\boldsymbol{\theta}}$.

### *3.2. NLP with ECs and IECs*

Now we consider the IECs in Problem 1 and establish the right DOE in $\mathbb{D}_f$, which not only satisfies the differential relation (22) but also meets the differential relation allowed by the IECs (2), i.e.,

$$\frac{\mathrm{d}g_i}{\mathrm{d}\tau} = (g_i)_{\boldsymbol{\theta}}{}^{\mathrm{T}} \frac{\mathrm{d}\boldsymbol{\theta}}{\mathrm{d}\tau} \le 0 \qquad i \in \mathbb{I} \ , \tag{41}$$

where $(g_i)_{\boldsymbol{\theta}}$ is the gradient of the $i$-th component of $\boldsymbol{g}$. $\mathbb{I}$ is the index set of the activated IECs for certain $\boldsymbol{\theta}$ and it is defined as

$$\mathbb{I} = \{ i \mid g_i(\boldsymbol{\theta}) = 0,\ i = 1,2,...,r \} \ . \tag{42}$$

In sub-section 3.1, we constructed the FPDOP to derive the DOE for the NLP with ECs only. Here the differential constraints (41) arising from the IECs (2) also need to be considered in constructing the FPDOP, which is now a typical QP problem as follows.

$$\begin{aligned} & \min \quad J_{t3} = \frac{1}{2} J_{t1} + \frac{1}{2} J_{t2} \\ & \text{subject to} \\ & \boldsymbol{h}_{\boldsymbol{\theta}} \frac{\mathrm{d}\boldsymbol{\theta}}{\mathrm{d}\tau} = \boldsymbol{0} \\ & (g_i)_{\boldsymbol{\theta}}{}^{\mathrm{T}} \frac{\mathrm{d}\boldsymbol{\theta}}{\mathrm{d}\tau} \le 0 \qquad i \in \mathbb{I} \end{aligned} \ , \tag{43}$$

where $J_{t1}$ and $J_{t2}$ are given by Eqs. (26) and (28), respectively. Of course we cannot simply set that $\frac{\mathrm{d}g_i}{\mathrm{d}\tau} = 0$ for all $i \in \mathbb{I}$ in the FPDOP, because such treatment may produce the wrong solution. From Definition 2.3, it is easy to find that strengthening an optimal active IEC to be an EC, the optimal solution will not be changed. Also, removing an IEC that is inactive at the optimal solution from the optimization problem, the optimal solution will not be changed either. Moreover, through the sensitivity theory [2], the multiplier information may also be deduced. See

**Lemma 3.1** [34]: *Strengthening an IEC to be an EC in the optimization problem, the corresponding multiplier is non-negative if this IEC is an optimal active IEC, and it is negative if this IEC is not an optimal active IEC.*

Therefore, to obtain the right DOE that may seek the optimal solution, in the FPDOP (43) only the optimal active IECs in (41) need to be considered. For the IECs in (41) that are not optimal active, the corresponding $g_i$ will fall into the inactive domain automatically. Introduce the index set of optimal active IECs for the FPDOP as follows

$$\mathbb{I}_p = \{ i \mid g_i = 0,\ \frac{\mathrm{d}g_i}{\mathrm{d}\tau} \le 0 \ \text{ is an optimal active IEC for the FPDOP (43)}, \quad i = 1,2,...,r \} \ , \tag{44}$$

with the number of its elements denoted by $n_{\mathbb{I}_p}$. Similarly, for the optimal solution $\hat{\boldsymbol{\theta}}$, $\mathbb{I}_p$ is highlighted with a hat "^" as $\hat{\mathbb{I}}_p$. We now present the following variant of the FPDOP (43) as

$$
\begin{aligned}
&\min \quad J_{t3}=\frac{1}{2}J_{t1}+\frac{1}{2}J_{t2} \\
&\text{subject to} \\
&\boldsymbol{h}_{\boldsymbol{\theta}}\frac{\mathrm{d}\boldsymbol{\theta}}{\mathrm{d}\tau}=\boldsymbol{0} \qquad\qquad . \\
&(g_i)_{\boldsymbol{\theta}}{}^{\mathrm{T}}\frac{\mathrm{d}\boldsymbol{\theta}}{\mathrm{d}\tau}=0 \qquad i\in\mathbb{I}_p
\end{aligned}
\tag{45}
$$

Through solving this problem analytically, we may obtain the DOE for general NLPs within $\mathbb{D}_f$ . See

**Theorem 3.3:** *For the NLP defined in Problem* 1*, under the premise of Assumption 2.1, solve the IVP defined by the following DOE*

$$
\frac{\mathrm{d}\boldsymbol{\theta}}{\mathrm{d}\tau}=-\boldsymbol{K}_{\boldsymbol{\theta}}(f_{\boldsymbol{\theta}}+\boldsymbol{h}_{\boldsymbol{\theta}}{}^{\mathrm{T}}\boldsymbol{\pi}_E+\boldsymbol{g}_{\boldsymbol{\theta}}{}^{\mathrm{T}}\boldsymbol{\pi}_I)\,,
\tag{46}
$$

*with arbitrary initial condition* $\boldsymbol{\theta}\big|_{\tau=0}\in\mathbb{D}_f$ *, where the parameter vectors* $\boldsymbol{\pi}_E\in\mathbb{R}^s$ *and* $\boldsymbol{\pi}_I\in\mathbb{R}^r$ *are determined by*

$$
\begin{aligned}
&\boldsymbol{\pi}_E=\begin{bmatrix}\pi_1\\ \pi_2\\ \cdots\\ \pi_s\end{bmatrix} \\
&(\pi_I)_i=0 \qquad i\notin\mathbb{I}_p,\;\; \boldsymbol{\pi}_I(\mathbb{I}_p)=\begin{bmatrix}\pi_{s+1}\\ \pi_{s+2}\\ \cdots\\ \pi_{s+n_{\mathbb{I}_p}}\end{bmatrix}\geq\boldsymbol{0}
\end{aligned}
\;,
\tag{47}
$$

*and the parameter vector* $\boldsymbol{\pi}\in\mathbb{R}^{s+n_{\mathbb{I}_p}}$ *is computed by*

$$
\boldsymbol{\pi}=-(\bar{\boldsymbol{h}}_{\boldsymbol{\theta}}\boldsymbol{K}_{\boldsymbol{\theta}}\bar{\boldsymbol{h}}_{\boldsymbol{\theta}}{}^{\mathrm{T}})^{+}\bar{\boldsymbol{h}}_{\boldsymbol{\theta}}\boldsymbol{K}_{\boldsymbol{\theta}}\boldsymbol{f}_{\boldsymbol{\theta}}\,.
\tag{48}
$$

$\boldsymbol{K}_{\boldsymbol{\theta}}$ *is an* $n\times n$ *dimensional positive-definite matrix and* $\bar{\boldsymbol{h}}=\begin{bmatrix}\boldsymbol{h}\\ \boldsymbol{g}(\mathbb{I}_p)\end{bmatrix}$. *When* $\tau\to+\infty$ , $\boldsymbol{\theta}$ *will satisfy the optimality condition* (7)*. Moreover, for the optimal solution,* $\hat{\mathbb{I}}_p=\hat{\mathbb{I}}$ *and the value of* $\boldsymbol{\pi}$ *is independent of* $\boldsymbol{K}_{\boldsymbol{\theta}}$ .

The proof is similar to Theorems 3.1 and 3.2. Regarding the argument that $\hat{\mathbb{I}}_p=\hat{\mathbb{I}}$ for the optimal solution of the NLP, this is because any component in $\hat{\mathbb{I}}$ also belongs to $\hat{\mathbb{I}}_p$ ultimately, or this active IEC will become inactive. Likewise, Theorem 3.3 implies the convergence of the multipliers, namely

**Remark 3.2:** Solve the NLP defined in Problem 1 with the approach suggested in Theorem 3.3; the Lagrange multiplier and the KKT multiplier given by Eqs. (47) and (48) will converge to the right solution as $\boldsymbol{\theta}$ converges to $\hat{\boldsymbol{\theta}}$ .

With the property of QP and since $\frac{\mathrm{d}\boldsymbol{\theta}}{\mathrm{d}\tau}=\boldsymbol{0}$ is always a feasible solution for the FPDOP, we may claim

**Remark 3.3:** The solution for the FPDOP (43) exists uniquely and it is the global minimum.

In pursuing the optimal solution under the dynamics governed by Eq. (46), the set $\mathbb{I}_p$ needs to be determined dynamically. Generally, $\mathbb{I}$ is easy to get. Thus we may first strengthen all IECs in $\mathbb{I}$ to get the corresponding Lagrange multipliers (and a feasible solution) and then use Lemma 3.1 to seek the right $\mathbb{I}_p$, following the procedure of the active-set methods [1-3].

The DOE proposed in this section is applicable when a feasible initial guess of the optimization parameter is available, and it may preserve the feasibility of solutions in the optimization. Therefore, the "merit function", which is usually employed for traditional NLP methods that cannot maintain the feasibility, is no longer necessary [3].

## 4. DOE valid in general

The DOE defined in Section 3 requires a feasible initial parameter vector in seeking the solution of the NLP. However, finding a feasible solution is usually not an easy task. In this section, we will generalize the DOE in the infeasible solution region $\mathbb{D}_{if}$ . The basic principle that we employ to eliminate the infeasibilities is the asymptotic stability of the first-order stable dynamic system; that is, an error parameter $e$ will be driven to zero in terms of the following equation

$$\frac{\mathrm{d}e}{\mathrm{d}\tau}=-ke\,, \tag{49}$$

where $k$ is a positive scalar. We will use Eq. (49) to address the problem of turning an infeasible solution that violates Eqs. (2) or (3) to be feasible.

### *4.1. Elimination of infeasibility on ECs*

Similarly, we start with a simple case that the IECs (2) are always satisfied strictly while the ECs (3) are violated, i.e.

$$\boldsymbol{g}(\boldsymbol{\theta})<\boldsymbol{0}\,, \tag{50}$$

$$\boldsymbol{h}(\boldsymbol{\theta})\neq\boldsymbol{0}\,. \tag{51}$$

Then we hope upon the DOE, the violated ECs (51) will achieve

$$\frac{\mathrm{d}\boldsymbol{h}}{\mathrm{d}\tau}=\boldsymbol{h}_{\boldsymbol{\theta}}\frac{\mathrm{d}\boldsymbol{\theta}}{\mathrm{d}\tau}=-\boldsymbol{K}_h\boldsymbol{h}\,, \tag{52}$$

where $\boldsymbol{K}_h$ is an $s\times s$ dimensional positive-definite matrix. Follow the strategy to get the DOE in Section 3; we may construct the Feasibility-achieving Dynamics Optimization Problem (FADOP) that may realize Eq. (52) as

$$\begin{aligned} &\min \quad J_{t3} = \frac{1}{2} J_{t1} + \frac{1}{2} J_{t2} \\ &\text{subject to} \\ &\frac{\mathrm{d}\boldsymbol{h}}{\mathrm{d}\tau} + \boldsymbol{K}_{\boldsymbol{h}} \boldsymbol{h} = \boldsymbol{0} \end{aligned} , \tag{53}$$

where $J_{t1}$ and $J_{t2}$ are given by Eqs. (26) and (28), respectively. Solve this FADOP analytically with the Lagrange multiplier technique; we may derive the DOE same to Eq. (23), while now the parameter vector $\boldsymbol{\pi}_E \in \mathbb{R}^s$, under the precondition that $\boldsymbol{h_\theta}$ has full row rank, is determined by

$$\boldsymbol{\pi}_E = -(\boldsymbol{h_\theta K_\theta h_\theta}^{\mathrm{T}})^{-1}(\boldsymbol{h_\theta K_\theta f_\theta} - \boldsymbol{K_h h}) . \tag{54}$$

If this precondition does not hold, not only the matrix $\boldsymbol{h_\theta K_\theta h_\theta}^{\mathrm{T}}$ is not invertible, but also no exact solution for Eq. (52) exists. Thus, certain modification is required.

**Proposition 4.1:** *With the parameter vector* $\boldsymbol{\pi}_E \in \mathbb{R}^s$ *determined by*

$$\boldsymbol{\pi}_E = -(\boldsymbol{h_\theta K_\theta h_\theta}^{\mathrm{T}})^{+}(\boldsymbol{h_\theta K_\theta f_\theta} - \boldsymbol{K_h h}) , \tag{55}$$

*the solution given by Eq.* (23) *is the optimal solution of the FADOP* (53) *with the constraint replaced with*

$$\frac{\mathrm{d}\boldsymbol{h}}{\mathrm{d}\tau} + \boldsymbol{P}_{\boldsymbol{h}_\theta c} \boldsymbol{K_h h} = \boldsymbol{0} , \tag{56}$$

*where* $\boldsymbol{P}_{\boldsymbol{h}_\theta c} = \boldsymbol{h_\theta h_\theta}^{+}$ *is the projection matrix.*

***Proof:*** Let $\boldsymbol{b} = \boldsymbol{P}_{\boldsymbol{h}_\theta c} \boldsymbol{K_h h}$, and obviously $\mathrm{rank}(\boldsymbol{h_\theta}) = \mathrm{rank}(\begin{bmatrix} \boldsymbol{h_\theta} & -\boldsymbol{b} \end{bmatrix})$. Remove the linearly dependent row vectors in $\boldsymbol{h_\theta}$ to obtain $\underline{\boldsymbol{h}}_{\boldsymbol{\theta}}$ and the same rows in $\boldsymbol{b}$ to obtain $\underline{\boldsymbol{b}}$. Then we get an equivalent expression of Eq. (56), in the sense that they have the same solution, as

$$\underline{\boldsymbol{h}}_{\boldsymbol{\theta}} \frac{\mathrm{d}\boldsymbol{\theta}}{\mathrm{d}\tau} = -\underline{\boldsymbol{b}} . \tag{57}$$

Replace the ECs in the FADOP (53) with Eq. (57) and denote the corresponding multiplier by $\underline{\boldsymbol{\pi}}_E$. With $\underline{\boldsymbol{\pi}}_E$ calculated by Eq. (54) and substituted into Eq. (23), the optimal solution to the reformulated FADOP is

$$\frac{\mathrm{d}\boldsymbol{\theta}}{\mathrm{d}\tau} = -\boldsymbol{K_\theta} \underline{\boldsymbol{h}}_{\boldsymbol{\theta}}^{\mathrm{T}} (\underline{\boldsymbol{h}}_{\boldsymbol{\theta}} \boldsymbol{K_\theta} \underline{\boldsymbol{h}}_{\boldsymbol{\theta}}^{\mathrm{T}})^{-1} \underline{\boldsymbol{b}} - (\boldsymbol{K_\theta} - \boldsymbol{K_\theta} \underline{\boldsymbol{h}}_{\boldsymbol{\theta}}^{\mathrm{T}} (\underline{\boldsymbol{h}}_{\boldsymbol{\theta}} \boldsymbol{K_\theta} \underline{\boldsymbol{h}}_{\boldsymbol{\theta}}^{\mathrm{T}})^{-1} \underline{\boldsymbol{h}}_{\boldsymbol{\theta}} \boldsymbol{K_\theta}) \boldsymbol{f_\theta} . \tag{58}$$

To prove this proposition, we only need to show that $\frac{\mathrm{d}\boldsymbol{\theta}}{\mathrm{d}\tau}$ given by Eq. (58) equals the following expression, i.e.,

$$\frac{\mathrm{d}\boldsymbol{\theta}}{\mathrm{d}\tau} = -\boldsymbol{K_\theta h_\theta}^{\mathrm{T}} (\boldsymbol{h_\theta K_\theta h_\theta}^{\mathrm{T}})^{+} \boldsymbol{K_h h} - (\boldsymbol{K_\theta} - \boldsymbol{K_\theta h_\theta}^{\mathrm{T}} (\boldsymbol{h_\theta K_\theta h_\theta}^{\mathrm{T}})^{+} \boldsymbol{h_\theta K_\theta}) \boldsymbol{f_\theta} , \tag{59}$$

which is obtained by substituting Eq. (55) into Eq. (23). The second terms in the right sides of Eqs. (58) and (59) are equal (see Eq. (39)). By Remark 2.3, we may further establish the equality of the first terms. See

$$\boldsymbol{K}_{\boldsymbol{\theta}}\underline{\boldsymbol{h}}_{\boldsymbol{\theta}}^{\mathrm{T}}(\underline{\boldsymbol{h}}_{\boldsymbol{\theta}}\boldsymbol{K}_{\boldsymbol{\theta}}\underline{\boldsymbol{h}}_{\boldsymbol{\theta}}^{\mathrm{T}})^{-1}\underline{\boldsymbol{b}}=\boldsymbol{K}_{\boldsymbol{\theta}}^{1/2}(\underline{\boldsymbol{h}}_{\boldsymbol{\theta}}\boldsymbol{K}_{\boldsymbol{\theta}}^{1/2})^{+}\underline{\boldsymbol{b}}=\boldsymbol{K}_{\boldsymbol{\theta}}^{1/2}(\boldsymbol{h}_{\boldsymbol{\theta}}\boldsymbol{K}_{\boldsymbol{\theta}}^{1/2})^{+}\boldsymbol{b}=\boldsymbol{K}_{\boldsymbol{\theta}}\boldsymbol{h}_{\boldsymbol{\theta}}^{\mathrm{T}}(\boldsymbol{h}_{\boldsymbol{\theta}}\boldsymbol{K}_{\boldsymbol{\theta}}\boldsymbol{h}_{\boldsymbol{\theta}}^{\mathrm{T}})^{+}\boldsymbol{b}=\boldsymbol{K}_{\boldsymbol{\theta}}\boldsymbol{h}_{\boldsymbol{\theta}}^{\mathrm{T}}(\boldsymbol{h}_{\boldsymbol{\theta}}\boldsymbol{K}_{\boldsymbol{\theta}}\boldsymbol{h}_{\boldsymbol{\theta}}^{\mathrm{T}})^{+}\boldsymbol{K}_{\boldsymbol{h}}\boldsymbol{h}\,, \tag{60}$$

where $\boldsymbol{K}_{\boldsymbol{\theta}}^{1/2}$ is the positive-definite matrix that satisfies $\boldsymbol{K}_{\boldsymbol{\theta}}^{1/2}\boldsymbol{K}_{\boldsymbol{\theta}}^{1/2}=\boldsymbol{K}_{\boldsymbol{\theta}}$. Thus, the proof is completed. □

Proposition 4.1 shows that with the pseudo-inverse, the error to the expected dynamics of the ECs may be minimized when $\boldsymbol{h}_{\boldsymbol{\theta}}$ do not have full row rank. Moreover, using the property of the projection matrix, it can be guaranteed that the ECs function values still tend to zero. This may be explained by

$$\frac{\mathrm{d}(\boldsymbol{h}^{\mathrm{T}}\boldsymbol{K}_{\boldsymbol{h}}\boldsymbol{h})}{\mathrm{d}\tau}=-2(\boldsymbol{K}_{\boldsymbol{h}}\boldsymbol{h})^{\mathrm{T}}\boldsymbol{P}_{\boldsymbol{h}_{\theta}c}\boldsymbol{K}_{\boldsymbol{h}}\boldsymbol{h}\leq 0\,. \tag{61}$$

Actually, Eq. (59) may be reformulated as

$$\frac{\mathrm{d}\boldsymbol{\theta}}{\mathrm{d}\tau}=-\boldsymbol{K}_{\boldsymbol{\theta}}^{1/2}(\boldsymbol{h}_{\boldsymbol{\theta}}\boldsymbol{K}_{\boldsymbol{\theta}}^{1/2})^{+}\boldsymbol{K}_{\boldsymbol{h}}\boldsymbol{h}-\boldsymbol{K}_{\boldsymbol{\theta}}^{1/2}\left(\boldsymbol{I}-(\boldsymbol{h}_{\boldsymbol{\theta}}\boldsymbol{K}_{\boldsymbol{\theta}}^{1/2})^{+}(\boldsymbol{h}_{\boldsymbol{\theta}}\boldsymbol{K}_{\boldsymbol{\theta}}^{1/2})\right)\boldsymbol{K}_{\boldsymbol{\theta}}^{1/2}\boldsymbol{f}_{\boldsymbol{\theta}}\,. \tag{62}$$

In order to understand Eq. (62) well, assuming $\boldsymbol{K}_{\boldsymbol{\theta}}=\boldsymbol{I}$ produces

$$\frac{\mathrm{d}\boldsymbol{\theta}}{\mathrm{d}\tau}=-\boldsymbol{h}_{\boldsymbol{\theta}}^{+}\boldsymbol{K}_{\boldsymbol{h}}\boldsymbol{h}-(\boldsymbol{I}-\boldsymbol{h}_{\boldsymbol{\theta}}^{+}\boldsymbol{h}_{\boldsymbol{\theta}})\boldsymbol{f}_{\boldsymbol{\theta}}\,. \tag{63}$$

Denote the row space spanned by the row vectors of $\boldsymbol{h}_{\boldsymbol{\theta}}$ with $\mathbb{S}_{\boldsymbol{h}_{\theta}r}$. Then the first term in the right hand of Eq. (63) indicates the solution projected to $\mathbb{S}_{\boldsymbol{h}_{\theta}r}$ from any solution of Eq. (56). Since $(\boldsymbol{I}-\boldsymbol{h}_{\boldsymbol{\theta}}^{+}\boldsymbol{h}_{\boldsymbol{\theta}})$ is the projection matrix to the orthogonal complement space (denoted by $\mathbb{S}_{\boldsymbol{h}_{\theta}r}^{\perp}$) of $\mathbb{S}_{\boldsymbol{h}_{\theta}r}$, the second term is the projection of the gradient $-\boldsymbol{f}_{\boldsymbol{\theta}}$ into $\mathbb{S}_{\boldsymbol{h}_{\theta}r}^{\perp}$, which is associated with the DOE (21) for the unconstrained optimization. For Eq. (62), it is just the generalized weighted case.

**Remark 4.1:** The DOE for the NLP with ECs includes two parts. One is the solution in $\mathbb{S}_{\boldsymbol{h}_{\theta}r}$ that eliminates the constraint violation; the other is the projection of objective function gradient in $\mathbb{S}_{\boldsymbol{h}_{\theta}r}^{\perp}$.

With the full row rank precondition of $\boldsymbol{h}_{\boldsymbol{\theta}}$, Yamashita [24] first established Eq. (54) and proved the convergence of the DOE solution to the optimal solution of the NLP. Here with Eq. (55), the ECs are expected to be satisfied gradually even if $\boldsymbol{h}_{\boldsymbol{\theta}}$ is not row full-rank. However, we will not prove the solution convergence until we resolve the violated IECs in a different way from former references.

*4.2. Elimination of infeasibility on IECs*

Now we expand the results to further accommodate the violated IECs, that is

$$\boldsymbol{g}(\boldsymbol{\theta}) > \boldsymbol{0} . \tag{64}$$

Analogously, we expect that for the violated IECs (and activated IECs), their dynamic motions satisfy

$$\frac{\mathrm{d} g_i}{\mathrm{d}\tau} = (g_i)_{\boldsymbol{\theta}}{}^{\mathrm{T}} \frac{\mathrm{d}\boldsymbol{\theta}}{\mathrm{d}\tau} \le -k_{gi} g_i \qquad i \in \mathbb{I} , \tag{65}$$

where the index set $\mathbb{I}$ is modified as

$$\mathbb{I} = \{ i \mid g_i(\boldsymbol{\theta}) \ge 0,\ i = 1,2,...,r \} , \tag{66}$$

and $k_{gi}$ is a positive constant for the $i$-th IEC. Likewise, the inequation (65) will be included in the FADOP as

$$\begin{aligned} &\min \quad J_{t3} = \frac{1}{2} J_{t1} + \frac{1}{2} J_{t2} \\ &\text{subject to} \\ &\quad \frac{\mathrm{d}\boldsymbol{h}}{\mathrm{d}\tau} + \boldsymbol{K}_{\boldsymbol{h}} \boldsymbol{h} = \boldsymbol{0} \\ &\quad \frac{\mathrm{d} g_i}{\mathrm{d}\tau} + k_{gi} g_i \le 0 \qquad i \in \mathbb{I} \end{aligned} \ . \tag{67}$$

With the precondition that the solution for the FADOP exists, we similarly introduce the modified index set $\mathbb{I}_p$ as

$$\mathbb{I}_p = \{ i \mid g_i \ge 0,\ \frac{\mathrm{d} g_i}{\mathrm{d}\tau} + k_{gi} g_i \le 0 \text{ is an optimal active IEC for the FADOP (67)}, i = 1,2,...,r \} , \tag{68}$$

and then the FADOP may be re-presented as

$$\begin{aligned} &\min \quad J_{t3} = \frac{1}{2} J_{t1} + \frac{1}{2} J_{t2} \\ &\text{subject to} \\ &\quad \frac{\mathrm{d}\boldsymbol{h}}{\mathrm{d}\tau} + \boldsymbol{K}_{\boldsymbol{h}} \boldsymbol{h} = \boldsymbol{0} \\ &\quad \frac{\mathrm{d}\boldsymbol{g}(\mathbb{I}_p)}{\mathrm{d}\tau} + \boldsymbol{K}_{\boldsymbol{g}} \boldsymbol{g}(\mathbb{I}_p) = \boldsymbol{0} \end{aligned} \ , \tag{69}$$

where $\boldsymbol{K}_{\boldsymbol{g}} = \mathrm{diag}\left(\boldsymbol{k}_{\boldsymbol{g}}(\mathbb{I}_p)\right)$ is an $n_{\mathbb{I}_p} \times n_{\mathbb{I}_p}$ dimensional positive-definite diagonal matrix and $\boldsymbol{k}_{\boldsymbol{g}} = \begin{bmatrix} k_{g1} & k_{g2} & \dots & k_{gr} \end{bmatrix}^{\mathrm{T}}$.

From the FADOP (69), we hope to obtain the DOE that may eliminate the violation on the IECs. However, for its rationality, one may argue that using the infinite-time asymptotic convergence principle, the violated IECs may never enter $\mathbb{D}_f$. Thus a finite-time convergence dynamics should be used in the FADOP such as

$$\frac{\mathrm{d}\boldsymbol{g}(\mathbb{I}_p)}{\mathrm{d}\tau} + \boldsymbol{K}_{\boldsymbol{g}} \mathrm{sign}\left(\boldsymbol{g}(\mathbb{I}_p)\right) = \boldsymbol{0} \ , \tag{70}$$

where sign(·) is the sign function. For this argument, we will show that even if an infinite-time convergence dynamics is employed for the IECs in $\mathbb{I}_P$, all the violated IECs in $\mathbb{I}$ will succeed in approaching the right region as long as the solution of the FADOP exists. Moreover, although Eq. (70) is usable theoretically, it may result in the unexpected chattering that is disadvantageous to the numerical computation.

**Proposition 4.2:** *Assuming its solution exists, the FADOP* (69) *guarantees that the infeasibilities on the IECs will be eliminated ultimately; concretely, the IECs with indexes in* $\hat{\mathbb{I}}$ *(See Eq.* (8) *for its definition) will at least achieve the feasibility asymptotically, and the IECs not in* $\hat{\mathbb{I}}$ *will return to the feasible solution region* $\mathbb{D}_f$ *in finite time.*

***Proof*****:** According to the definition of the FADOP (69), for the unviolated IECs, they will always stay in the feasible solution region $\mathbb{D}_f$. Therefore, we only consider the violated IECs in $\mathbb{I}$. For optimal active IECs of the FADOP (69), there are

$$\frac{\mathrm{d}g_i}{\mathrm{d}\tau}+k_{gi}g_i=0 \qquad i\in\mathbb{I}_p, \tag{71}$$

and for inactive IECs of the FADOP, we have

$$\frac{\mathrm{d}g_i}{\mathrm{d}\tau}<-k_{gi}g_i \qquad i\notin\mathbb{I}_p. \tag{72}$$

On the other hand, an activated IEC $g_i$, which is not optimal active, will enter the inactive domain from $g_i=0$ with

$$\frac{\mathrm{d}g_i}{\mathrm{d}\tau}<-a \ \ when \ \ g_i=0, \tag{73}$$

where $a$ is some positive constant. Thus, there exists a neighbourhood $[0,\varepsilon]$ such that

$$\frac{\mathrm{d}g_i}{\mathrm{d}\tau}<-\frac{a}{2}<-k_{gi}g_i \quad when \ \ g_i\in[0,\varepsilon], \tag{74}$$

where $\varepsilon$ is a small positive constant. Under the dynamics of Eq. (71) or Eq. (72), a violated IEC $g_i$ will enter $[0,\varepsilon]$ in finite time from any positive value. After that, an IEC $g_i$ that is not optimal active will reach zero with time smaller than $\frac{2\varepsilon}{a}$ from $g_i=\varepsilon$. Thus, for those IECs whose indexes are not in $\hat{\mathbb{I}}$ (i.e., $\hat{\mathbb{I}}_p$ as well), they will return to $\mathbb{D}_f$ in a limited time. For the IECs with indexes in $\hat{\mathbb{I}}$, if their indexes always belong to $\mathbb{I}_p$, their feasibility will be achieved asymptotically; if not, they may enter $\mathbb{D}_f$ in finite time. □

Through solving the FADOP (69), we may derive the DOE same to Eq. (46), with the parameters $\boldsymbol{\pi}_E$ and $\boldsymbol{\pi}_I$ also determined by Eq. (47) but $\boldsymbol{\pi}$ computed by

$$\boldsymbol{\pi} = -(\bar{\boldsymbol{h}}_{\boldsymbol{\theta}} \boldsymbol{K}_{\boldsymbol{\theta}} \bar{\boldsymbol{h}}_{\boldsymbol{\theta}}^{\mathrm{T}})^{+} (\bar{\boldsymbol{h}}_{\boldsymbol{\theta}} \boldsymbol{K}_{\boldsymbol{\theta}} \boldsymbol{f}_{\boldsymbol{\theta}} - \begin{bmatrix} \boldsymbol{K}_h \boldsymbol{h} \\ \boldsymbol{K}_g \boldsymbol{g}(\mathbb{I}_p) \end{bmatrix}) , \tag{75}$$

where $\bar{\boldsymbol{h}} = \begin{bmatrix} \boldsymbol{h} \\ \boldsymbol{g}(\mathbb{I}_p) \end{bmatrix}$. Note that the pseudo-inverse computation may address the linear-dependence in rows of $\bar{\boldsymbol{h}}_{\boldsymbol{\theta}}$ and eliminate the infeasibilities arising from the desired dynamics for the violated constraints. With the resulting DOE, $\dfrac{\mathrm{d}\boldsymbol{\theta}}{\mathrm{d}\tau}$ is explicitly comprised of two parts as Remark 4.1 indicated. Besides, it implicitly ensures that the violations on the IECs, whose index is in $\mathbb{I}$ but not in $\mathbb{I}_p$, will be relieved as well. In particular, when $\bar{\boldsymbol{h}}_{\boldsymbol{\theta}}$ has full column rank, $\boldsymbol{I} - (\bar{\boldsymbol{h}}_{\boldsymbol{\theta}} \boldsymbol{K}_{\boldsymbol{\theta}}^{1/2})^{+} (\bar{\boldsymbol{h}}_{\boldsymbol{\theta}} \boldsymbol{K}_{\boldsymbol{\theta}}^{1/2}) = \boldsymbol{0}$, and then the DOE is equivalent to the following expression, i.e.,

$$\frac{\mathrm{d}\boldsymbol{\theta}}{\mathrm{d}\tau} = -\boldsymbol{K}_{\boldsymbol{\theta}}^{1/2} (\bar{\boldsymbol{h}}_{\boldsymbol{\theta}} \boldsymbol{K}_{\boldsymbol{\theta}}^{1/2})^{+} \begin{bmatrix} \boldsymbol{K}_h \boldsymbol{h} \\ \boldsymbol{K}_g \boldsymbol{g}(\mathbb{I}_p) \end{bmatrix}, \tag{76}$$

which minimizes the error to the expected dynamics for the violated constraints regardless of $\boldsymbol{f}_{\boldsymbol{\theta}}$. Again, the set $\mathbb{I}_p$ needs to be determined dynamically, and Lemma 3.1 may be used to seek the right results, as the active-set methods do. In particular, in seeking $\mathbb{I}_p$, the linearly dependent rows in $\boldsymbol{h}_{\boldsymbol{\theta}}$ need to be removed temporarily.

### *4.3. Mathematic validation*

With the modification in the last section, it is anticipated that the generalized DOE, defined by Eqs. (46), (47) and (75), will evolve an arbitrary initial guess of $\boldsymbol{\theta}$ to the optimal solution, achieving the feasibility and the optimality simultaneously. It is not hard to verify that the solution satisfying the feasibility conditions (2), (3) and the optimality condition (7) is the equilibrium solution of the generalized DOE. However, when considered in the infeasible solution region $\mathbb{D}_{if}$, the objective (1) cannot act as the Lyapunov function to ensure the convergence of the solution. See; the change of the objective function value, i.e., Eq. (20), may be re-presented as

$$\frac{\mathrm{d}J}{\mathrm{d}\tau} = (\boldsymbol{f}_{\boldsymbol{\theta}} + \boldsymbol{h}_{\boldsymbol{\theta}}^{\mathrm{T}} \boldsymbol{\pi}_E + \boldsymbol{g}_{\boldsymbol{\theta}}^{\mathrm{T}} \boldsymbol{\pi}_I)^{\mathrm{T}} \frac{\mathrm{d}\boldsymbol{\theta}}{\mathrm{d}\tau} - \boldsymbol{\pi}_E^{\mathrm{T}} \frac{\mathrm{d}\boldsymbol{h}}{\mathrm{d}\tau} - \boldsymbol{\pi}_I^{\mathrm{T}} \frac{\mathrm{d}\boldsymbol{g}}{\mathrm{d}\tau} . \tag{77}$$

By investigating Eq. (77), it is found that starting from an infeasible solution, $J$ will not monotonically decrease under Eqs. (46), (47) and (75), because the term $-\boldsymbol{\pi}_E^{\mathrm{T}} \dfrac{\mathrm{d}\boldsymbol{h}}{\mathrm{d}\tau} - \boldsymbol{\pi}_I^{\mathrm{T}} \dfrac{\mathrm{d}\boldsymbol{g}}{\mathrm{d}\tau}$ arising from the infeasibilities may be positive, and the sign of $\dfrac{\mathrm{d}J}{\mathrm{d}\tau}$ is uncertain.

Lacking the theoretical guarantee, it is natural to ask that is it ensured that $\boldsymbol{\theta}$ will approach the equilibrium solution from arbitrary initial value, rather than converge to the bounded limit cycle as the Van der Pol oscillator (in which the equilibrium is unstable, see Khalil [38]). Now we will answer this question with rigorous mathematic argument, and a Lyapunov function will be constructed. Before we carry out the mathematic analysis, certain assumptions are presented.

**Assumption 4.1:** The solution for the FADOP exists.

**Assumption 4.2:** The multiplier parameters $\boldsymbol{\pi}_E$ and $\boldsymbol{\pi}_I$, determined by Eqs. (47) and (75), are bounded as

$$\left\|\boldsymbol{\pi}_E\right\|_2 \le d_E , \tag{78}$$

$$\left\|\boldsymbol{\pi}_I\right\|_2 \le d_I , \tag{79}$$

where $d_E$ and $d_I$ are some positive constants.

**Lemma 4.1:** *For the function*

$$V = \sqrt{\boldsymbol{h}^{\mathrm{T}}\boldsymbol{h}} + \sqrt{\boldsymbol{g}(\mathbb{I})^{\mathrm{T}}\boldsymbol{g}(\mathbb{I})} + c_1 J , \tag{80}$$

*where $J$ is defined in Eq.* (1) *and $\mathbb{I}$ is defined in Eq.* (66)*, on the basis of Assumptions 2.1, 4.1 and 4.2, there exists certain positive constant $c_1$*

$$c_1 < \min\left(\frac{\lambda_{\min}(\boldsymbol{K_h})}{d_E \lambda_{\max}(\boldsymbol{K_h})}, \frac{\min(\boldsymbol{k_g})}{d_I \max(\boldsymbol{k_g})}\right), \tag{81}$$

*such that Eq.* (80) *is a Lyapunov function for the generalized DOE defined by Eqs.* (46), (47) *and* (75). *Here $\lambda_{\min}(\cdot)$ and $\lambda_{\max}(\cdot)$ denote the minimum eigenvalue and the maximum eigenvalue of matrix, respectively.*

***Proof:*** First, we show that the minimum solution of Problem 1 is also the minimum solution of the function (80), and vice versa. If the parameters $\boldsymbol{\theta}$ are located within the feasible solution region $\mathbb{D}_f$, then we have

$$V = c_1 J . \tag{82}$$

Obviously for this case, the minimum of Problem 1 and the minimum of the function (80) are the same. When the parameters lie in the infeasible solution region $\mathbb{D}_{if}$, we consider the neighbourhood around the minimum solution $\hat{\boldsymbol{\theta}}$. Since $\hat{\boldsymbol{\theta}}$ satisfies Eq. (7), we have the first-order expansion of the function (80) around $\hat{\boldsymbol{\theta}}$ as

$$\mathrm{d}V = \left\|\mathrm{d}\boldsymbol{h}\right\|_2 + \left\|\mathrm{d}\boldsymbol{g}(\mathbb{I})\right\|_2 - c_1 \boldsymbol{\pi}_E{}^{\mathrm{T}} \mathrm{d}\boldsymbol{h} - c_1 \boldsymbol{\pi}_I(\mathbb{I}_p)^{\mathrm{T}} \mathrm{d}\boldsymbol{g}(\mathbb{I}_p) . \tag{83}$$

Here note that $\mathrm{d}J$ is obtained from Eq. (77) and

$$\boldsymbol{\pi}_I{}^{\mathrm{T}}\mathrm{d}\boldsymbol{g} = \boldsymbol{\pi}_I(\mathbb{I}_p)^{\mathrm{T}}\mathrm{d}\boldsymbol{g}(\mathbb{I}_p) \ . \tag{84}$$

According to Assumption 4.2, and with the Holder's inequality, there are

$$-d_E\left\|\mathrm{d}\boldsymbol{h}\right\|_2 \le -\left\|\boldsymbol{\pi}_E\right\|_2\left\|\mathrm{d}\boldsymbol{h}\right\|_2 \le -\boldsymbol{\pi}_E{}^{\mathrm{T}}\mathrm{d}\boldsymbol{h} \ , \tag{85}$$

$$-d_I\left\|\mathrm{d}\boldsymbol{g}(\mathbb{I})\right\|_2 \le -\left\|\boldsymbol{\pi}_I\right\|_2\left\|\mathrm{d}\boldsymbol{g}(\mathbb{I}_p)\right\|_2 \le -\boldsymbol{\pi}_I(\mathbb{I}_p)^{\mathrm{T}}\mathrm{d}\boldsymbol{g}(\mathbb{I}_p) \ . \tag{86}$$

Then we have

$$\mathrm{d}V \ge (1-c_1 d_E)\left\|\mathrm{d}\boldsymbol{h}\right\|_2 + (1-c_1 d_I)\left\|\mathrm{d}\boldsymbol{g}(\mathbb{I})\right\|_2 \ . \tag{87}$$

According to Eq. (81), we have $\mathrm{d}V > 0$. Thus, the solution $\hat{\boldsymbol{\theta}}$ determines a minimum for the function (80). On the other hand, the minimum solution of the function (80) is also the minimum solution of Problem 1 under Assumption 4.1, because if any infeasible solution for Problem 1 introduces a minimum for the function (80), there will be no solution exists for the FADOP to get rid of the violated constraints.

Now we consider the derivative of $V$ with respect to the virtual time $\tau$. Differentiating Eq. (80) with respect to $\tau$ produces

$$\frac{\mathrm{d}V}{\mathrm{d}\tau} = \frac{\boldsymbol{h}^{\mathrm{T}}}{\sqrt{\boldsymbol{h}^{\mathrm{T}}\boldsymbol{h}}}\frac{\mathrm{d}\boldsymbol{h}}{\mathrm{d}\tau} + \frac{\boldsymbol{g}(\mathbb{I})^{\mathrm{T}}}{\sqrt{\boldsymbol{g}(\mathbb{I})^{\mathrm{T}}\boldsymbol{g}(\mathbb{I})}}\frac{\mathrm{d}\boldsymbol{g}(\mathbb{I})}{\mathrm{d}\tau} + c_1\frac{\mathrm{d}J}{\mathrm{d}\tau} \ . \tag{88}$$

Substitute Eq. (77) in and use Eqs. (46), (47), and (75), and especially note that

$$\boldsymbol{g}(\mathbb{I})^{\mathrm{T}}\frac{\mathrm{d}\boldsymbol{g}(\mathbb{I})}{\mathrm{d}\tau} \le -\boldsymbol{g}(\mathbb{I})^{\mathrm{T}}\mathrm{diag}\left(\boldsymbol{k}_{\boldsymbol{g}}(\mathbb{I})\right)\boldsymbol{g}(\mathbb{I}) \ , \tag{89}$$

$$\boldsymbol{\pi}_I{}^{\mathrm{T}}\frac{\mathrm{d}\boldsymbol{g}}{\mathrm{d}\tau} = \boldsymbol{\pi}_I(\mathbb{I}_p)^{\mathrm{T}}\boldsymbol{K}_{\boldsymbol{g}}\boldsymbol{g}(\mathbb{I}_p) \ . \tag{90}$$

According to Assumption 4.1, we have

$$\begin{aligned}\frac{\mathrm{d}V}{\mathrm{d}\tau} \le & -\frac{\boldsymbol{h}^{\mathrm{T}}}{\sqrt{\boldsymbol{h}^{\mathrm{T}}\boldsymbol{h}}}\boldsymbol{K}_{\boldsymbol{h}}\boldsymbol{h} - \frac{\boldsymbol{g}(\mathbb{I})^{\mathrm{T}}}{\sqrt{\boldsymbol{g}(\mathbb{I})^{\mathrm{T}}\boldsymbol{g}(\mathbb{I})}}\mathrm{diag}\left(\boldsymbol{k}_{\boldsymbol{g}}(\mathbb{I})\right)\boldsymbol{g}(\mathbb{I}) - c_1(f_{\boldsymbol{\theta}} + \boldsymbol{h}_{\boldsymbol{\theta}}{}^{\mathrm{T}}\boldsymbol{\pi}_E + \boldsymbol{g}_{\boldsymbol{\theta}}{}^{\mathrm{T}}\boldsymbol{\pi}_I)^{\mathrm{T}}\boldsymbol{K}(f_{\boldsymbol{\theta}} + \boldsymbol{h}_{\boldsymbol{\theta}}{}^{\mathrm{T}}\boldsymbol{\pi}_E + \boldsymbol{g}_{\boldsymbol{\theta}}{}^{\mathrm{T}}\boldsymbol{\pi}_I) \\ & + c_1\boldsymbol{\pi}_E{}^{\mathrm{T}}\boldsymbol{K}_{\boldsymbol{h}}\boldsymbol{h} + c_1\boldsymbol{\pi}_I(\mathbb{I}_p)^{\mathrm{T}}\boldsymbol{K}_{\boldsymbol{g}}\boldsymbol{g}(\mathbb{I}_p)\end{aligned} \ . \tag{91}$$

According to Assumption 4.2, and with the Holder's inequality, we have

$$\boldsymbol{\pi}_E{}^{\mathrm{T}}\boldsymbol{K}_{\boldsymbol{h}}\boldsymbol{h} \le \lambda_{\max}(\boldsymbol{K}_{\boldsymbol{h}})\left\|\boldsymbol{\pi}_E\right\|_2\left\|\boldsymbol{h}\right\|_2 \le d_E\lambda_{\max}(\boldsymbol{K}_{\boldsymbol{h}})\left\|\boldsymbol{h}\right\|_2 \ . \tag{92}$$

In particular, there is

$$\boldsymbol{\pi}_I(\mathbb{I}_p)^{\mathrm{T}}\boldsymbol{K}_{\boldsymbol{g}}\boldsymbol{g}(\mathbb{I}_p) \le \lambda_{\max}(\boldsymbol{K}_{\boldsymbol{g}})\left\|\boldsymbol{\pi}_I\right\|_2\left\|\boldsymbol{g}(\mathbb{I}_p)\right\|_2 \le d_I\max(\boldsymbol{k}_{\boldsymbol{g}})\left\|\boldsymbol{g}(\mathbb{I})\right\|_2 \ . \tag{93}$$

Substituting the inequalities (92) and (93) into Eq. (91) gives

$$
\begin{aligned}
\frac{\mathrm{d}V}{\mathrm{d}\tau} \le & -\left(\lambda_{\min}(\boldsymbol{K}_{\boldsymbol{h}}) - c_1 d_E \lambda_{\max}(\boldsymbol{K}_{\boldsymbol{h}})\right)\|\boldsymbol{h}\|_2 - \left(\min(\boldsymbol{k}_{\boldsymbol{g}}) - c_1 d_I \max(\boldsymbol{k}_{\boldsymbol{g}})\right)\|\boldsymbol{g}(\mathbb{I})\|_2 \\
& - c_1 (f_{\boldsymbol{\theta}} + {\boldsymbol{h}_{\boldsymbol{\theta}}}^{\mathrm{T}}\boldsymbol{\pi}_E + {\boldsymbol{g}_{\boldsymbol{\theta}}}^{\mathrm{T}}\boldsymbol{\pi}_I)^{\mathrm{T}} \boldsymbol{K} (f_{\boldsymbol{\theta}} + {\boldsymbol{h}_{\boldsymbol{\theta}}}^{\mathrm{T}}\boldsymbol{\pi}_E + {\boldsymbol{g}_{\boldsymbol{\theta}}}^{\mathrm{T}}\boldsymbol{\pi}_I)
\end{aligned} . \tag{94}
$$

With $c_1$ set by Eq. (81), we have that $\frac{\mathrm{d}V}{\mathrm{d}\tau} \le 0$ hold in both $\mathbb{D}_f$ and $\mathbb{D}_{if}$, and $\frac{\mathrm{d}V}{\mathrm{d}\tau} = 0$ when Eqs. (2), (3), and (7) are satisfied. □

**Theorem 4.1:** *For the NLP defined in Problem 1, solve the IVP defined by the DOE* (46) *with arbitrary initial condition* $\boldsymbol{\theta}\big|_{\tau=0}$*, where the parameters* $\boldsymbol{\pi}_E \in \mathbb{R}^s$ *and* $\boldsymbol{\pi}_I \in \mathbb{R}^r$ *are determined by Eqs.* (47) *and* (75)*. Then under Assumptions 2.1, 4.1 and 4.2, when* $\tau \to +\infty$*,* $\boldsymbol{\theta}$ *will satisfy the feasibility conditions* (2)*,* (3) *and the optimality condition* (7)*.*

***Proof:*** The proof is a direct application of Lemma 2.2 and Lemma 4.1. From Lemma 4.1, the function (80) is ensured a Lyapunov function for the dynamic system, defined by Eqs. (46), (47), and (75), around the equilibrium that meets (2), (3), and (7). According to Lemma 2.2, the equilibrium solution is an asymptotically stable solution. Therefore $\boldsymbol{\theta}$ will satisfy the feasibility conditions (2), (3) and the optimality condition (7) of the NLP defined in Problem 1 when $\tau \to +\infty$. □

In particular, Theorem 4.1 means a global convergence property of the DOE, since the claim that the function (80) is a valid Lyapunov function holds for arbitrarily large region around the minimums with $c_1$ small enough.

**Remark 4.2:** Presume Assumptions 2.1, 4.1 and 4.2 hold. Then solve the NLP defined in Problem 1 with the approach suggested in Theorem 4.1; the solution will converge to the optimal solution globally.

Similarly, we may claim the convergence of the multipliers, which is implied by Theorem 4.1.

**Remark 4.3:** Presume Assumptions 2.1, 4.1 and 4.2 hold. Then solve the NLP defined in Problem 1 with the approach suggested in Theorem 4.1; the Lagrange multiplier and the KKT multiplier given by Eqs. (47) and (75) will converge to the right solution as $\boldsymbol{\theta}$ converges to $\hat{\boldsymbol{\theta}}$.

The generalized DOE, developed in this section, may operate in both the feasible and the infeasible solution region. In the feasible solution region, it may better preserve the feasibility than the DOE in Section 3 because error arising from the numerical computation will be eliminated. When it is applied in the infeasible solution region, the Lyapunov function (80) may play the role of a "merit function", which measures the improvement of solutions. However, this merit function only has theoretical significance and will not be used in the computation.

### *4.4. Discussion*

Active-set methods to solve the FADOP require a feasible initial solution. When the number of the constraints is not larger than $n$, generating a feasible solution is easy, simply by equalizing all IECs. When there are many violated IECs, an auxiliary Linear Programming (LP) problem may be set up to seek the feasible solution for the FADOP, that is

$$\begin{aligned} & \min \quad \gamma \\ & \text{subject to} \\ & \boldsymbol{h}_\theta \frac{\mathrm{d}\boldsymbol{\theta}}{\mathrm{d}\tau} + \boldsymbol{K}_h \boldsymbol{h} = \boldsymbol{0} \\ & \frac{(g_i)_\theta^{\mathrm{T}}}{\|(g_i)_\theta\|_2} \frac{\mathrm{d}\boldsymbol{\theta}}{\mathrm{d}\tau} + k_{gi} \frac{g_i}{\|(g_i)_\theta\|_2} \leq \gamma \qquad i \in \mathbb{I} \end{aligned} \quad . \tag{95}$$

The proposed LP has clear geometric meaning, in which the parameter $\gamma$ represents the minimum distance (with sign) to the IECs. Solutions for the LP always exist and the unbounded case may be avoided by restricting the range of $\frac{\mathrm{d}\boldsymbol{\theta}}{\mathrm{d}\tau}$ or $\gamma$. The solution $\frac{\mathrm{d}\boldsymbol{\theta}}{\mathrm{d}\tau}$ from the LP (95) is a feasible solution for the FADOP (67) if $\gamma \leq 0$, and no feasible solution for the FADOP (67) exists if $\gamma > 0$. Since a feasible initial guess for this LP is easy to obtain, its solution may be solved efficiently with the simplex method.

A pivotal premise that we establish the global convergence of the generalized DOE is Assumption 4.1. However, in practice this assumption may be violated. For the infeasibility arising from the ECs, it may be addressed with the projection matrix. For the infeasibilities coming from the IECs, the parameter $\gamma$ of LP (95) actually provides the reference margin for the tuning of $k_{gi}$. However, the violation on Assumption 4.1 may still occur and then the solution of the DOE will halt at the wrong results. To ensure that the developed DOE works well under arbitrary initial conditions, technique of Priority Treatment Strategy (PTS) on the IECs may be employed. In the computation, the PTS treats the IECs according to their priority. Part of the IECs will be considered first and the others will be temporarily ignored. After achieving the feasibilities of these prior IECs, the rest will be included orderly, and all the IECs will be met in a sequential way. With the PTS, Assumption 4.1 is easier to be met and the optimal solution will be sought effectively.

## 5. Examples

First, an unconstrained problem is studied to investigate the performance of the elementary gradient-based dynamic method; that is, the optimization on the Generalized Wood function (GENWOOD) [41].

**Example 5.1:**

$$\min \ f(\boldsymbol{\theta}) = 1 + \sum_{i\in\Lambda}\begin{pmatrix} 100(\theta_{i+1}-\theta_i^2)^2 + (1-\theta_i)^2 + 90(\theta_{i+3}-\theta_{i+2}{}^2)^2 \\ +(1-\theta_{i+2})^2 + 10(\theta_{i+1}+\theta_{i+3}-2)^2 + 0.1(\theta_{i+1}-\theta_{i+3})^2 \end{pmatrix}$$

where $\boldsymbol{\theta}\in\mathbb{R}^n$ is the optimization parameter vector, $n$ is a multiple of 4 and $\Lambda=\{1,5,9,...,n-3\}$. The optimal solution for this example is $\hat{\boldsymbol{\theta}}=\begin{bmatrix}1 & 1 & ... & 1\end{bmatrix}_{n\times 1}{}^{\mathrm{T}}$.

With the dynamic method, the resulting DOE is given by Eq. (21), where $f_{\boldsymbol{\theta}}=\begin{bmatrix}\boldsymbol{d}_1{}^{\mathrm{T}} & \boldsymbol{d}_2{}^{\mathrm{T}} & ... & \boldsymbol{d}_j{}^{\mathrm{T}} & ... & \boldsymbol{d}_{n/4}{}^{\mathrm{T}}\end{bmatrix}_{n\times 1}{}^{\mathrm{T}}$ and

$$\boldsymbol{d}_j=\begin{bmatrix} -400(\theta_{i+1}-\theta_i^2)\theta_i-2(1-\theta_i) \\ 200(\theta_{i+1}-\theta_i^2)+20(\theta_{i+1}+\theta_{i+3}-2)+0.2(\theta_{i+1}-\theta_{i+3}) \\ -360(\theta_{i+3}-\theta_{i+2}{}^2)\theta_{i+2}-2(1-\theta_{i+2}) \\ 180(\theta_{i+3}-\theta_{i+2}{}^2)+20(\theta_{i+1}+\theta_{i+3}-2)-0.2(\theta_{i+1}-\theta_{i+3}) \end{bmatrix}; j=1,2,...,\frac{n}{4}; i=4j-3$$

$\boldsymbol{K}_{\boldsymbol{\theta}}$ was set as $\boldsymbol{K}_{\boldsymbol{\theta}}=\boldsymbol{I}_{n\times n}$. Consider the case of $n=8$. The initial conditions of $\boldsymbol{\theta}$ were given by $\boldsymbol{\theta}\big|_{\tau=0}=\begin{bmatrix}-3 & -1 & -3 & -1 & -2 & 0 & -2 & 0\end{bmatrix}^{\mathrm{T}}$. Because the dynamics of the DOE are fast and a little stiff, to enhance the numerical stability when solving the resulting IVPs, the stiff ODE integrator "ode15s" in MATLAB, with a relative error tolerance of $1\times10^{-3}$ and an absolute error tolerance of $1\times10^{-6}$, is employed for the numerical integration along the time horizon of $[0,100]\mathrm{s}$, which might be extended adaptively.

Figure 2 plots the profiles of the parameters, showing that the solutions converge to the optimal after a heavy oscillation and the minimizer point is obtained after $\tau=60\mathrm{s}$. In particular, during the time span of $\tau\in[2,15]$ s, the gradient of the objective function is nearly zero and the corresponding parameters are almost unchanged around $\boldsymbol{\theta}\big|_{\tau=2.9}$ = [-0.9692 0.9494 -0.9683 0.9490 -0.9693 0.9498 -0.9681 0.9486]$^{\mathrm{T}}$. Furthermore, comparison with common iterative methods for unconstrained problems were made, including the Broyden-Fletcher-Goldfarb-Shanno (BFGS) quasi-Newton method, the trust-region method and the polytope method, by invoking the functions of "fminuc" and "fminsearch" in MATLAB. The results are summarized in Table 1. It is found that the BFGS quasi-Newton method and the polytope method failed to seek the optimal solution, but they both terminated successfully with tolerance criterion met. In particular, the quasi-Newton method halts at the results that are similar to the solutions around $\tau=2.9\mathrm{s}$ from the dynamic method. The trust-region method succeeds to get the optimal solution. However, the time consumption is larger than the dynamic method.

We further considered the case of $n=100$, which brings much more optimization parameters. The initial conditions of $\boldsymbol{\theta}$ were $\boldsymbol{\theta}\big|_{\tau=0}$ = [-3 -1 -3 -1 -2 0 -2 0 … … -2 0]$^{\mathrm{T}}$. The computation of the dynamic method finished rapidly with only 0.2674s consumed. The solution curves are presented in Fig. 3 and the overall profiles are very similar to those in Fig. 2. In contrast, the BFGS quasi-Newton method and the polytope method again failed while the trust-region method succeeds. However, the time it used is 1.0023s, much larger than that of the dynamic method, which shows the superior efficiency of the dynamic method.

Table 1. Computation results for Example 1.

| | Optimized parameters | Time consumption | Exit information |
|---|---|---|---|
| The dynamic method | [1.0000 1.0000 1.0000 1.0000 1.0000 1.0000 1.0000 1.0000]$^{T}$ | 0.0960s | The solution converges with parameter tolerance of $1\times10^{-4}$ under the given integral time horizon. |
| The BFGS quasi-Newton method | [-0.9556 0.9233 -0.9817 0.9750 -0.9691 0.9494 -0.9684 0.9490]$^{T}$ | 0.0168s | Magnitude of gradient is smaller than the optimality tolerance of $1\times10^{-6}$. |
| The trust-region method | [1.0000 1.0000 1.0000 1.0000 1.0000 1.0000 1.0000 1.0000]$^{T}$ | 0.1907s | Change in the objective function value was less than the function tolerance of $1\times10^{-6}$. |
| The polytope method | [0.9226 0.8507 1.0770 1.1603 0.2143 0.0280 0.1680 -0.0010]$^{T}$ | 0.0412s | The function converged to a solution with function tolerance of $1\times10^{-4}$ and parameter tolerance of $1\times10^{-4}$. |

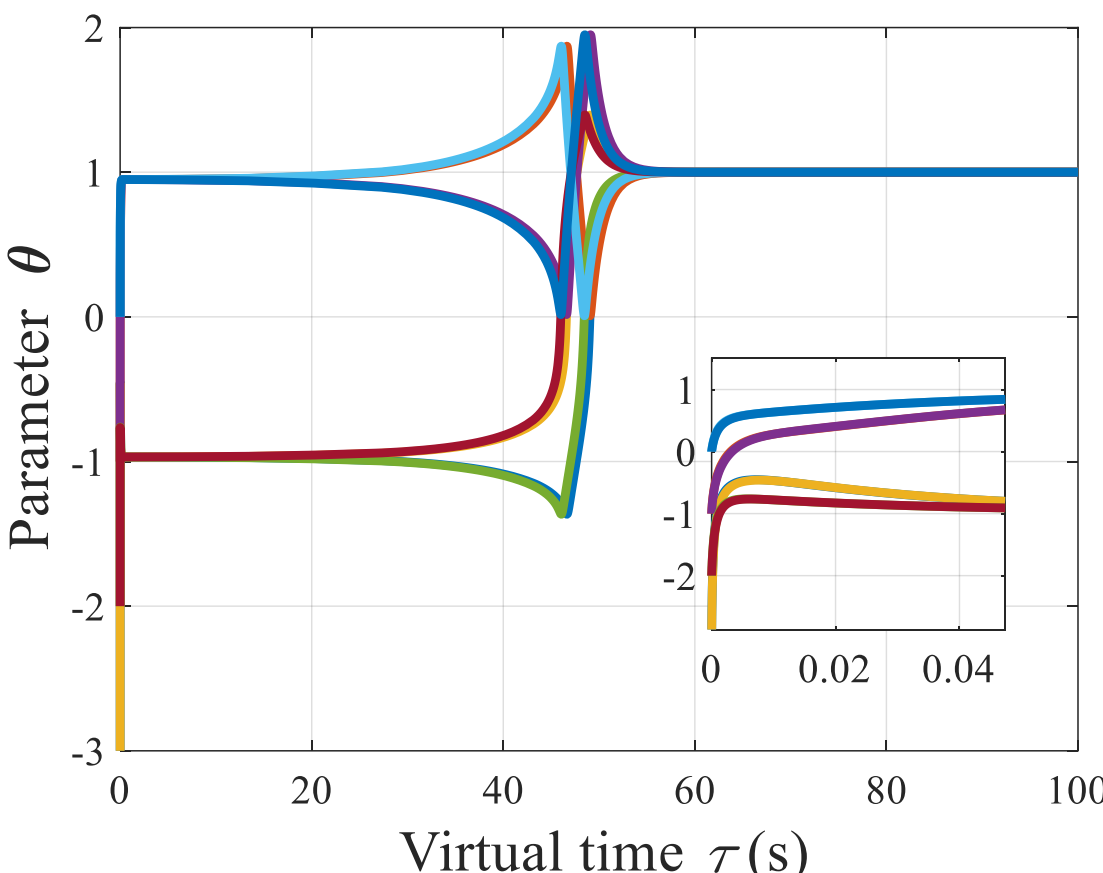

Figure 2. The dynamic motion curves of the optimization parameters for Example 1 ($n$=8).

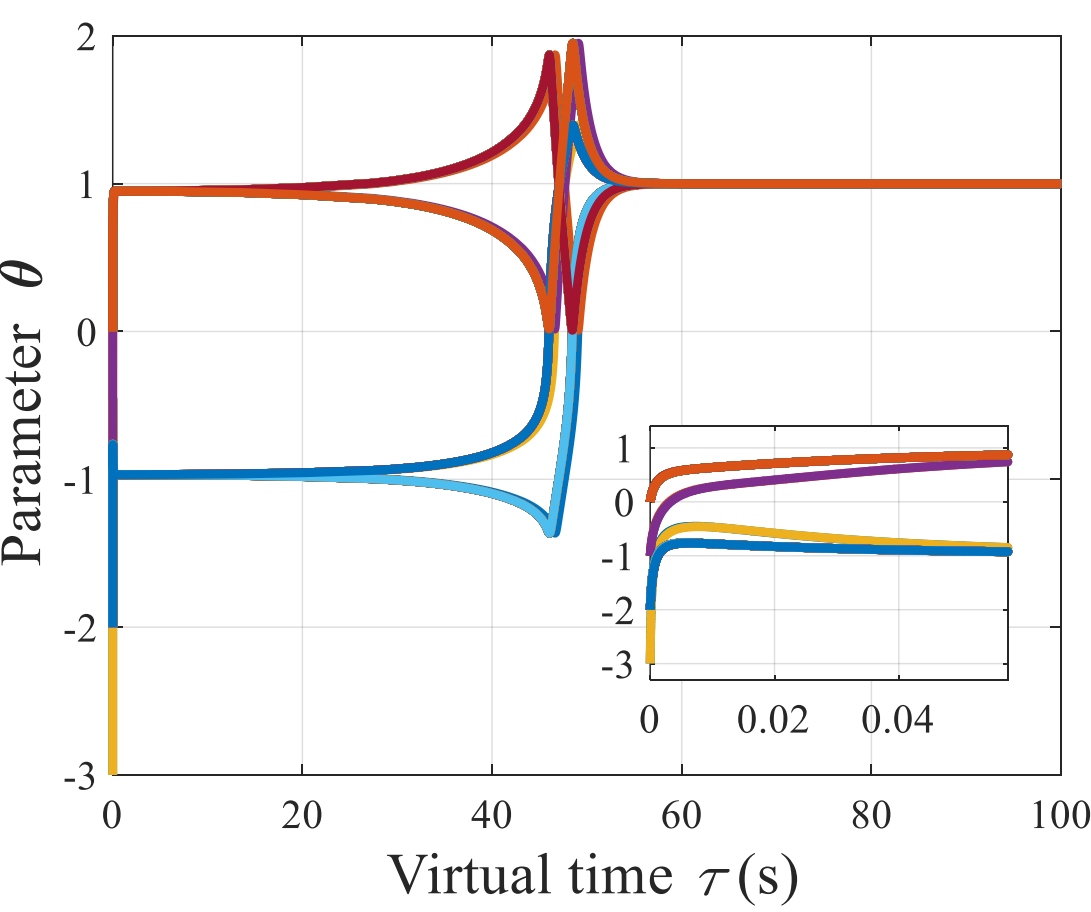

Figure 3. The dynamic motion curves of the optimization parameters for Example 1 ($n$=100).

Next a NLP adapted from Li et al. [42] is considered and note that a redundant EC is intentionally included.

**Example 5.2:**

$$\min \quad f(\boldsymbol{\theta}) = -\theta_1\theta_2 - \theta_2\theta_3 - \theta_3\theta_1$$

subject to

$$g_1(\boldsymbol{\theta}) = -\theta_1 \le 0 \qquad g_2(\boldsymbol{\theta}) = -\theta_2 \le 0 \qquad g_3(\boldsymbol{\theta}) = -\theta_3 \le 0$$

$$g_4(\boldsymbol{\theta}) = 0.5(\theta_1 - 3)^2 + \theta_2^2 + \theta_3^2 - 1 \le 0 \qquad g_5(\boldsymbol{\theta}) = \frac{\theta_1}{0.5 + \theta_2^2} + 2\theta_3 - 4 \le 0$$

$$h_1(\boldsymbol{\theta}) = \theta_1 + \theta_2 + \theta_3 - 3 = 0 \qquad h_2(\boldsymbol{\theta}) = 2\theta_1 + 2\theta_2 + 2\theta_3 - 6 = 0$$

Using the proposed dynamic method, the terms in the DOE (46) are

$$\boldsymbol{f_\theta} = \begin{bmatrix} -\theta_2 - \theta_3 \\ -\theta_1 - \theta_3 \\ -\theta_2 - \theta_1 \end{bmatrix}, \boldsymbol{h_\theta} = \begin{bmatrix} 1 & 1 & 1 \\ 2 & 2 & 2 \end{bmatrix}, \boldsymbol{g_\theta} = \begin{bmatrix} -1 & 0 & 0 \\ 0 & -1 & 0 \\ 0 & 0 & -1 \\ \theta_1 - 3 & 2\theta_2 & 2\theta_3 \\ \dfrac{1}{0.5+\theta_2^2} & \dfrac{-2\theta_1\theta_2}{(0.5+\theta_2^2)^2} & 2 \end{bmatrix}.$$

The three-dimensional matrix $\boldsymbol{K_\theta}$ was set as $\boldsymbol{K_\theta} = 0.1 \cdot \boldsymbol{I}_{3\times3}$. The matrix $\boldsymbol{K_h}$ was $\boldsymbol{K_h} = 0.1 \cdot \boldsymbol{I}_{2\times2}$ and the vector parameter $\boldsymbol{k_g}$ was $\boldsymbol{k_g} = \begin{bmatrix} 0.1 & 0.1 & 0.1 & 0.1 & 0.1 \end{bmatrix}^{\mathrm{T}}$. PTS was used with the IECs $g_i$ $(i = 1,2,3)$ given higher priority. We tried 10 computations starting from different initial conditions of $\boldsymbol{\theta}|_{\tau=0}$ and their values were sampled from a uniform distribution within [-10, 10]. To solve the resulting IVPs, the ODE integrator "ode45" in MATLAB, with a relative error tolerance of $1\times10^{-3}$ and an absolute error tolerance of $1\times10^{-6}$, was employed on an integration time horizon of 300s.

Table 2 gives the computation results, in which $e_{\boldsymbol{\theta}} = \left\| \boldsymbol{\theta}|_{\tau=300} - \hat{\boldsymbol{\theta}} \right\|$ and the optimal solution is $\hat{\boldsymbol{\theta}} = \begin{bmatrix} 2 & 0.5 & 0.5 \end{bmatrix}^{\mathrm{T}}$. It is shown that the precision of the solution is high and the time used for the computation is small. For one case of $\boldsymbol{\theta}|_{\tau=0} = \begin{bmatrix} -4.8578 & 3.8180 & -2.7364 \end{bmatrix}^{\mathrm{T}}$, figure 4 gives the profiles of $\boldsymbol{\theta}(\tau)$. It is shown that they approach the optimal solutions quickly. Figure 5 gives the profiles of the multipliers for the ECs and the IECs $g_i$ $(i = 4,5)$. At $\tau = 300$ s, we compute that $(\pi_E)_1|_{\tau=300} = 0.35$, $(\pi_E)_2|_{\tau=300} = 0.70$, $(\pi_I)_4|_{\tau=300} = 0.75$, and $(\pi_I)_5|_{\tau=300} = 0$. They are consistent with the results in Ref. [42]. For the IEC $g_5$, the value of its multiplier is always zero, while for the IEC $g_4$, its multiplier gradually approaches the right solution during the dynamic process.

Table 2. Computation results for Example 1.

| | Average | Minimum | Maximum |
|---|---|---|---|
| Error $e_{\boldsymbol{\theta}}$ | $2.7395\times10^{-12}$ | $3.3171\times10^{-13}$ | $5.2464\times10^{-12}$ |
| Time consumed (s) | 0.0503 | 0.0437 | 0.0602 |

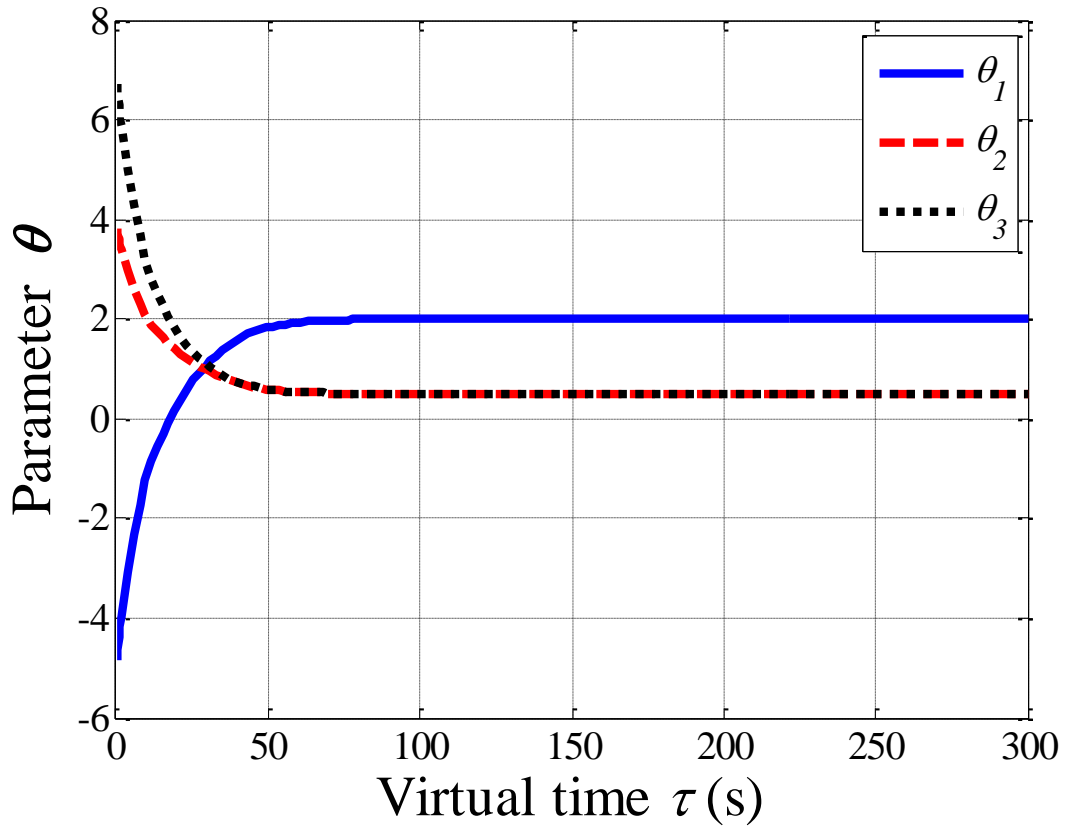

Figure 4. The dynamic motion curves of the optimization parameters for Example 1.

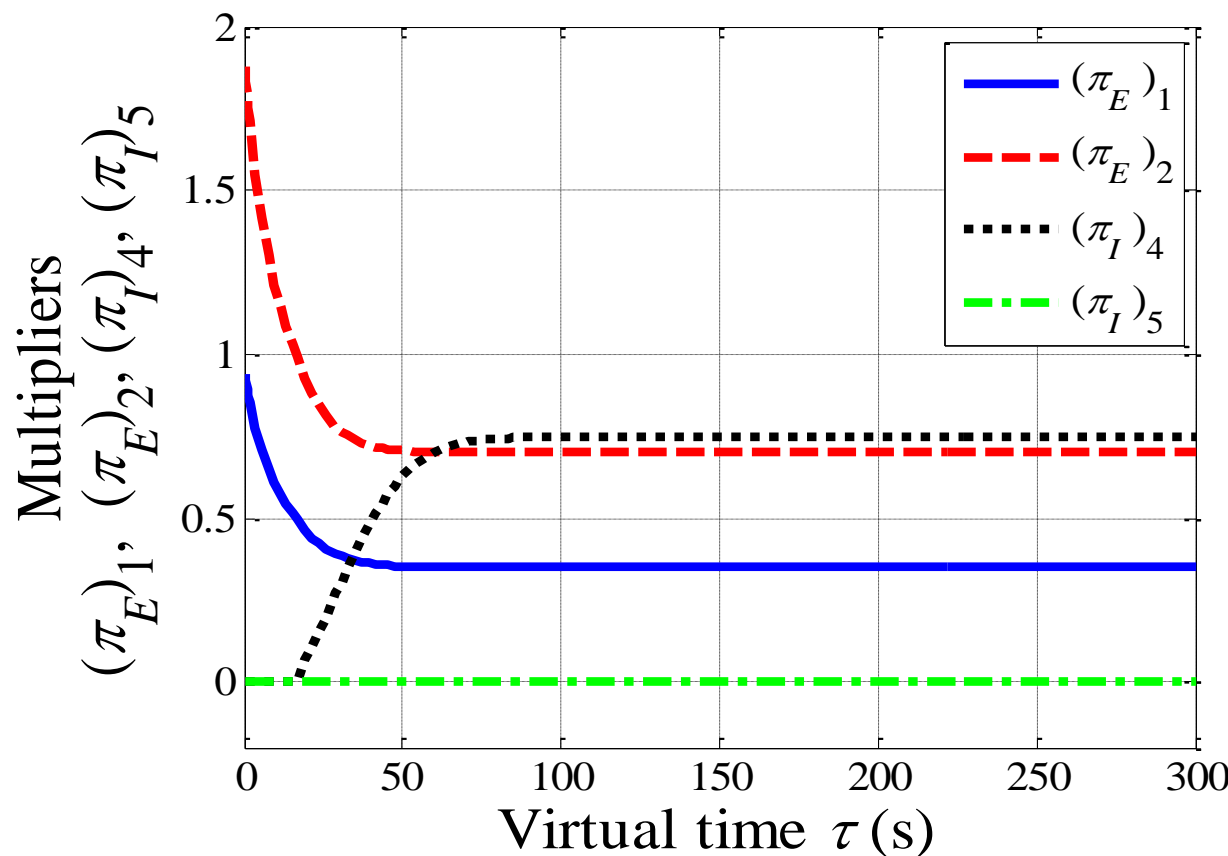

Figure 5. The profiles of the Lagrange multipliers $(\pi_E)_1$, $(\pi_E)_2$ and the KKT multipliers $(\pi_I)_4$, $(\pi_I)_5$ for Example 1.

Now we consider another NLP problem with much more optimization parameters and constraints from Dong [43].

**Example 5.3:**

$$\min \quad f(\boldsymbol{\theta}) = \sin(\theta_1 - 1 + 1.5\pi) + \sum_{2 \le i \le 100} 100\sin(-\theta_i + 1.5\pi + \theta_{i+1}^2)$$

subject to

$$\begin{aligned} 0.5 \le \theta_i \le 1.5 \qquad & i = 1 \\ -\pi \le \theta_{i-1}^2 - \theta_i \le \pi \qquad & i = 2,3,\ldots,100 \\ \theta_i - \theta_{i+1} = 0 \qquad & i = 1,2,\ldots,99 \end{aligned}$$

To solve this example, the parameter matrixes in the DOE were set as $\boldsymbol{K}_{\boldsymbol{\theta}} = 0.1 \cdot \boldsymbol{I}_{100\times100}$ and $\boldsymbol{K}_{\boldsymbol{h}} = 1 \cdot \boldsymbol{I}_{99\times99}$. The IECs were rearranged to the standard form and the vector parameter $\boldsymbol{k}_{\boldsymbol{g}}$ was $\boldsymbol{k}_{\boldsymbol{g}} = \left(\begin{bmatrix} 1 & 1 & \ldots & 1 \end{bmatrix}_{1\times200}\right)^{\mathrm{T}}$. We also tried 10 computations starting from different initial conditions, in which $\theta_1\big|_{\tau=0} = 2$, intentionally set out of the feasible region, and $\theta_i\big|_{\tau=0}$ $(i = 2,3,\ldots,100)$ were sampled from a uniform distribution within [0.7, 1.2]. Since there are many parameters in the DOE and this equation has a multiple-time-scale structure, the stiff ODE integrator

"ode15s" in MATLAB is employed for the numerical integration along the time horizon of $[0,100]$s . In the integrator setting, the relative error tolerance and absolute error tolerance were $1\times10^{-3}$ and $1\times10^{-6}$, respectively.

In Table 3, relevant computation results are presented. Here the optimal optimization parameters are $\hat{\boldsymbol{\theta}}=\left(\begin{bmatrix}1 & 1 & \dots & 1\end{bmatrix}_{1\times100}\right)^{\mathrm{T}}$ . It is again found that the accurate numerical solutions are obtained rapidly. Especially, for one designed case of $\theta_i\big|_{\tau=0}=2-\dfrac{3(i-1)}{99}$ $(i=1,2,\dots,100)$ , figure 6 gives the motion curves of $\boldsymbol{\theta}$ in the dynamic process. It is shown that they approach the value of 0.5 at first and then "jump" to the optimal value of 1. We further investigate the multipliers for the constraints. Figure 7 gives the profiles of the multipliers for the EC of $\theta_1-\theta_2=0$ and the IEC of $\theta_1-1.5\le0$ . It is found that the Lagrange multiplier oscillates heavily at the beginning and then converges to zero rapidly, while the KKT multiplier (even if this IEC is violated at $\tau=0$ s) always maintains the value of zero during the entire dynamic optimization process.

Table 3. Computation results for Example 2.

| | Average | Minimum | Maximum |
|---|---|---|---|
| Error $e_{\boldsymbol{\theta}}$ | $1.1138\times10^{-5}$ | $9.1476\times10^{-10}$ | $4.5180\times10^{-5}$ |
| Time consumed (s) | 1.2000 | 1.1712 | 1.2321 |

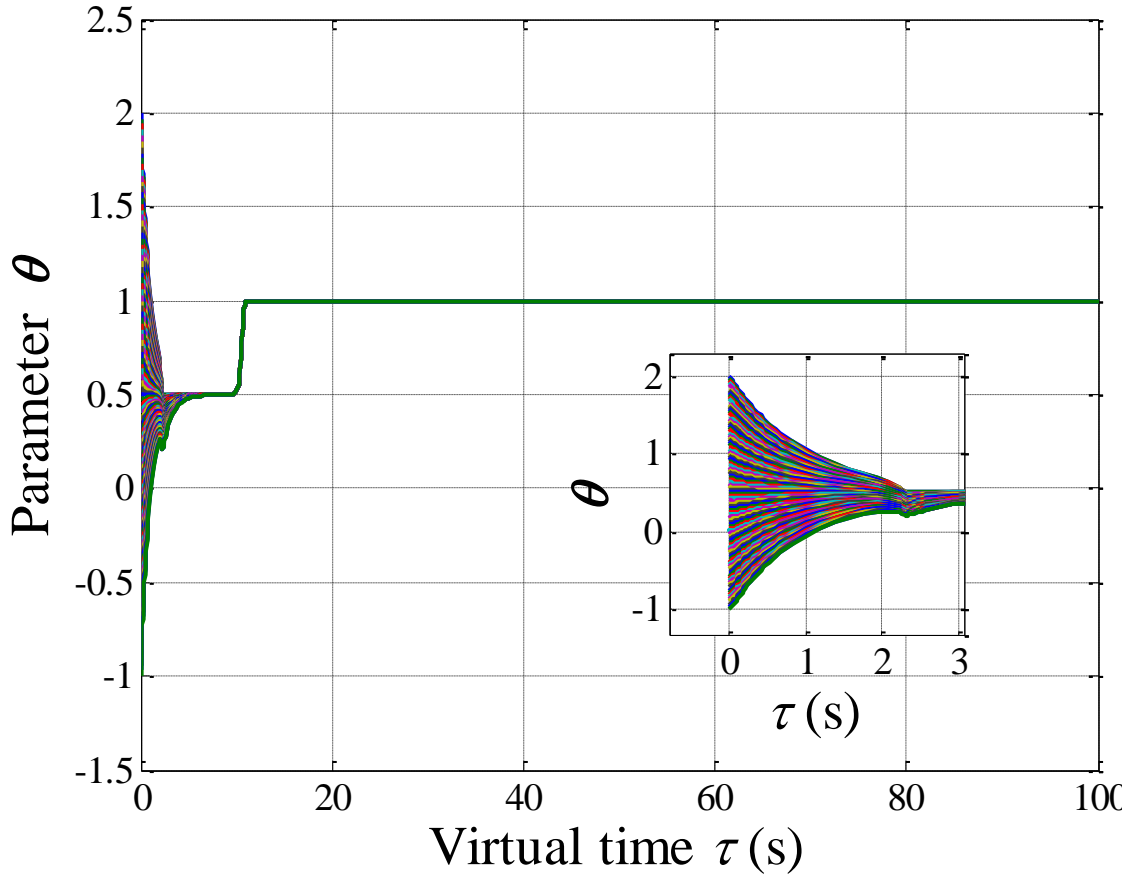

Figure 6. The dynamic motion curves of the optimization parameters for Example 2. There are 100 curves corresponding to $\theta_i$ $(i=1,2,\dots,100)$ .

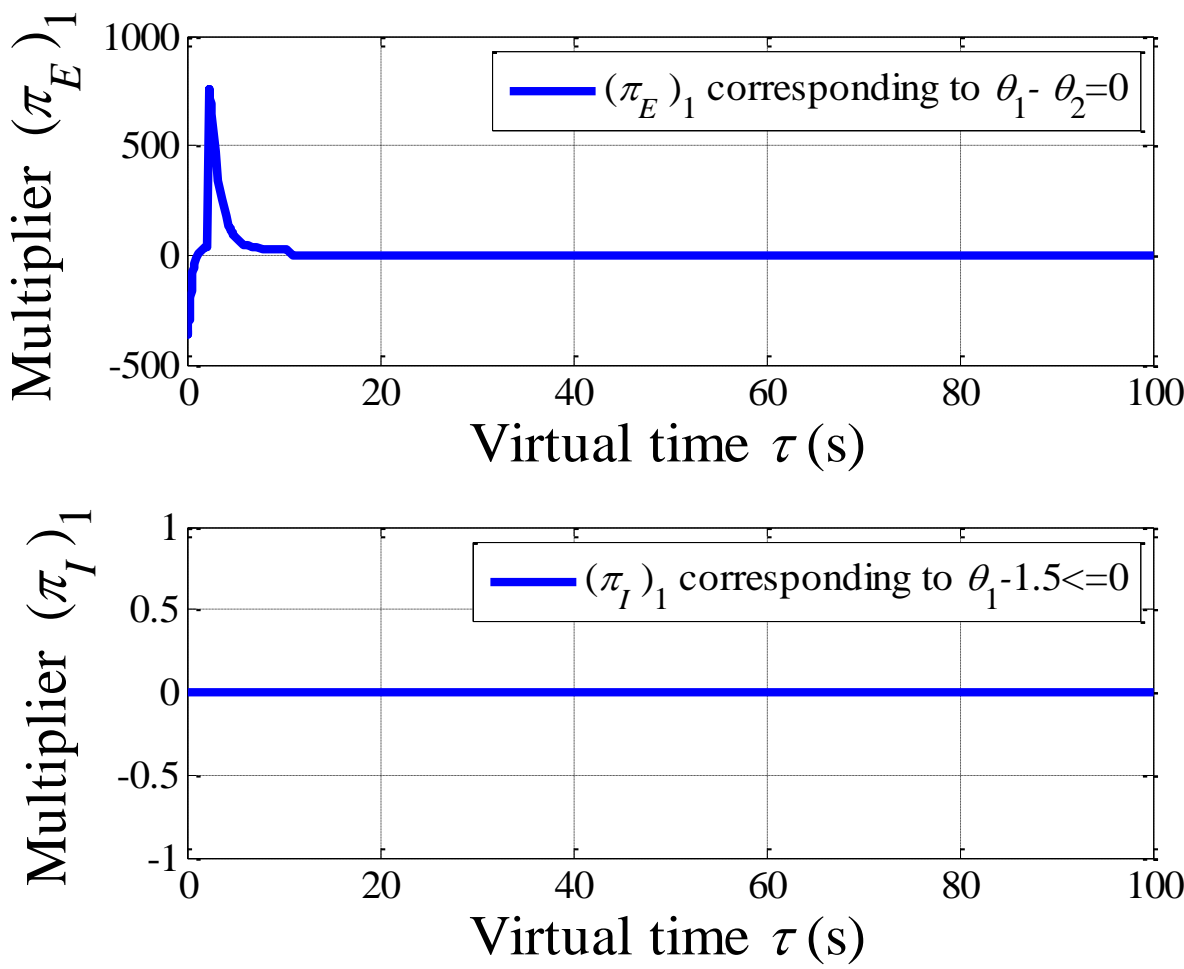

Figure 7. The profiles of the Lagrange multiplier $(\pi_E)_1$ and the KKT multiplier $(\pi_I)_1$ for Example 2.

## 6. Concluding remarks

Via the proposed dynamic method, the optimization of Non-linear Programming (NLP) problems is transformed to the computation of Initial-value Problems (IVPs) with arbitrary initial conditions. Upon the Lyapunov stability theory, the solution of the Dynamic Optimization Equation (DOE) is guaranteed to meet the optimality condition with global convergence, under the premise that the Karush-Kuhn-Tucker (KKT) condition exists. Note that although only the first-order optimality conditions are explicitly guaranteed, the solution will not halt at a maximum or a saddle (unless they are the initial guesses), since those solutions are not asymptotically stable equilibriums. Particularly with the pseudo-inverse, the DOEs proposed are also valid even if the linearly independent regularity condition on the nonlinear constraints is violated. Illustrative examples are solved, and it is shown that this method performs well in both the precision and the efficiency.

One may sense that the dynamic method and the traditional numerical iterative method are not completely irrelevant, in that our method shares some similarity to the Sequential Quadratic Programming (SQP) method as they both solve the Quadratic Programming (QP) sub-problem. Actually Schropp has already proved in his work that SQP method for the NLP with the Equality Constraints (ECs) can be regarded as a variable step size Euler-Cauchy integration method applied to the DOE. Intuitively, any discrete iterative method may have its continuous dynamic counterpart. In contrast to the iterative methods, with the dynamic method, the daunting task of searching for a reasonable step size and the annoying oscillation phenomenon around the optimum are eliminated, and the mature Ordinary Differential Equation (ODE) integration methods may be employed to solve the resulting IVP conveniently. However, the dynamic method is far less popular than the numerical iterative method. In our opinion, this is maybe because the iterative methods have become efficient and conventional in solving the NLP problems.

Yet remarkably, the dynamic method has a concise theoretical expression in the DOE form, and it may be more appropriate for the application in the control community as the augmented internal dynamics.